\documentclass[10pt]{amsart}
\usepackage{latexsym}   
\usepackage{amssymb}    
\usepackage{amsmath}    
\usepackage{amsthm}     
\usepackage[all]{xy}
\usepackage{hyperref}
\newcommand{\pa}{\partial}\newcommand{\al}{\alpha}

\newcommand{\del}{\delta}

\newcommand{\om}{\omega}

\newcommand{\ti}{\tilde}

\renewcommand{\thefootnote}

\newtheorem{theorem}{Theorem}[section]

\theoremstyle{definition}
\newtheorem{definition}[theorem]{Definition}

\theoremstyle{remark}

\numberwithin{equation}{section}

\title[On a new definition of the B\"{a}cklund transformation in the isometric deformation of surfaces]
{On a new definition of the B\"{a}cklund transformation in the isometric deformation of surfaces}

\author[  Ion I. Dinc\u{a}]{Ion I. Dinc\u{a}}
\address{Department of Applied Mathematics, Faculty of Applied Sciences,
National University of Science and Technology Politehnica Bucharest 313 Spl. Independentei 060042
Bucharest, Romania}
 \email{idinca@upb.ro}

\begin{document}


\begin{abstract}
We prove that a generic $4$-dimensional integrable rolling distribution of contact elements with the
symmetry of the tangency configuration (excluding developable seed and isotropic developable leaves) splits
into an $1$-dimensional family of generic $3$-dimensional integrable rolling distributions of contact elements
with the symmetry of the tangency configuration, thus introducing a new definition of the B\"{a}cklund transformation
in the isometric deformation of surfaces.
\end{abstract}

\maketitle

\tableofcontents \pagenumbering{arabic}

\section{Introduction}

The classical problem of finding the isometric deformations of surfaces (see Eisenhart \cite{E1})
was stated in 1859 by the French Academy of Sciences as

\begin{center}
{\it To find all surfaces applicable to a given one.}
\end{center}

Probably the most successful researcher of this problem is Bianchi, who in 1906 in \cite{B1}
solved the problem for quadrics by introducing the B\"{a}cklund transformation of surfaces
isometric to quadrics and the isometric correspondence provided by the Ivory affine transformation.
By 1909 Bianchi had a fairly complete treatment \cite{B2} and in 1917 he proved in \cite{B3} the
rigidity of the B\"{a}cklund transformation of isometric deformations of quadrics in the case of
auxiliary surface plane: the only B\"{a}cklund transformation with defining surface and auxiliary
surface plane appears as the singular B\"{a}cklund transformation of isometric deformations of
quadrics.

In \cite{D3} we have managed to improve Bianchi's result \cite{B3} to arbitrary auxiliary surface:
the only B\"{a}cklund transformation with defining surface is Bianchi's B\"{a}cklund transformation
of isometric deformations of quadrics.

We shall consider the complexification
$$(\mathbb{C}^3,\langle\cdot,\cdot\rangle),\ \langle x,y\rangle :=x^Ty,\
|x|^2:=x^Tx,\ x,y\in\mathbb{C}^3$$
of the real $3$-dimensional Euclidean space; in this setting surfaces are $2$-dimensional objects of
$\mathbb{C}^3$ depending on two real or complex parameters.

{\it Isotropic} (null) vectors are those vectors of length $0$; since most vectors are not isotropic
we call a vector simply vector and we shall emphasize isotropic for isotropic vectors. The same
denomination will apply in other settings: for example we call quadric a non-degenerate quadric (a
quadric projectively equivalent to the complex unit sphere).

Consider Lie's viewpoint: one can replace a surface $x\subset\mathbb{C}^3$ with a $2$-dimensional
distribution of {\it contact elements} (pairs of points and planes passing through those points; the
classical geometers call them {\it facets}): the collection of its tangent planes (with the points
of tangency highlighted); thus a contact element is the infinitesimal version of a surface (the
integral element $(x,dx)|_{\mathrm{pt}}$ of the surface). Conversely, a $2$-dimensional distribution
of contact elements is not always the collection of the tangent planes of a surface (with the points
of tangency highlighted), but the condition that a $2$-dimensional distribution of contact elements
is integrable (that is it is the collection of the tangent planes of a leaf (sub-manifold)) does not
distinguish between the cases when this sub-manifold is a surface, curve or point, thus allowing the
collapsing of the leaf.

A $3$-dimensional distribution of contact elements is integrable if it is the collection of the
tangent planes of an $1$-dimensional family of leaves.

Two {\it rollable} (applicable or isometric) surfaces can be {\it rolled} (applied) one onto the
other such that at any instant they meet tangentially and with same differential at the tangency
point.

\begin{definition}
The rolling of two isometric surfaces $x_0,\ x\subset\mathbb{C}^3$ (that is $|dx_0|^2=|dx|^2$) is
the surface, curve or point $(R,t)\subset\mathbf{O}_3(\mathbb{C})\ltimes\mathbb{C}^3$ such that
$(x,dx)=(R,t)(x_0,dx_0):=(Rx_0+t,Rdx_0)$.
\end{definition}

The rolling introduces the flat connection form (it encodes the difference of the second fundamental
forms of $x_0,\ x$ and it being flat encodes the difference of the Gau\ss-Codazzi-Mainardi-
Peterson equations of $x_0,\ x$).

\begin{definition}
Consider an integrable $3$-dimensional distribution of contact elements $\mathcal{F}=(p,m)$ centered
at $p=p(u,v,w)$, with normal fields $m=m(u,v,w)$ and distributed along the surface $x_0=x_0(u,v)$.
If we roll $x_0$ on an isometric surface $x$ (that is $(x,dx)=(R,t)(x_0,dx_0):=(Rx_0+t,Rdx_0)$), then
the rolled distribution of contact elements is $(Rp+t,Rm)$ and is distributed along $x$; if it
remains integrable for any rolling, then the distribution is called {\it integrable rolling
distribution of contact elements}.
\end{definition}

Since by infinitesimal rolling in an arbitrary tangential infinitesimal direction $\del$ an initial
contact element $\mathcal{F}$ which is common tangent plane to two isometric surfaces is replaced
with an infinitesimally close contact element $\mathcal{F}'$ having in common with $\mathcal{F}$ the
direction $\del$, in the actual rolling problem we have contact elements centered on each other (the
symmetric {\it tangency configuration}) and contact elements centered on another one $\mathcal{F}$
reflect in $\mathcal{F}$; note that we assumed a finite law of a general nature as a consequence of
an infinitesimal law via discretization (the converse is clear); see Bianchi \cite{B3}.

Thus for a theory of isometric deformation of surfaces with the assumption above we are led to
consider, via rolling, certain $4$-dimensional distributions of contact elements centered on the
tangent planes of the considered surface $x_0$ and passing through the origin of the tangent planes
(each point of each tangent plane is the center of finitely many contact elements) and their rolling
counterparts on the isometric surface $x$.

After a thorough study of infinitesimal laws (in particular infinitesimal isometric deformations)
and their iterations the classical geometers were led to consider the B\"{a}cklund transformation (a
finite law of a general nature) as a consequence of an infinitesimal law via discretization.

The B\"{a}cklund transformation in the isometric deformation problem naturally appears by splitting
the $4$-dimensional distribution of contact elements above into an $1$-dimensional family of
$3$-dimensional integrable rolling distributions of contact elements, thus introducing a spectral
parameter $z$ (the B\"{a}cklund transformation is denoted $B_z$); each $3$-dimensional distribution
of contact elements is integrable (with {\it leaves} $x^1$) regardless of the shape of the
{\it seed} surface $x^0$.

The question appears if for the $4$-dimensional distribution of contact elements above being
integrable for any surface $x$ isometric to $x_0$ and with a $2$-dimensional family of leaves (in
general surfaces), then it splits into an $1$-dimensional family of $3$-dimensional integrable
rolling distributions of contact elements.

The purpose of this note is to provide an affirmative answer to this question.

Bianchi considered the most general form of a B\"{a}cklund transformation as the focal surfaces (one
transform of the other) of a {\it Weingarten} congruence (congruence upon whose two focal surfaces
the asymptotic directions correspond; equivalently the second fundamental forms are proportional).
Note that although the correspondence provided by the Weingarten congruence does not give the
applicability (isometric) correspondence, the B\"{a}cklund transformation is the tool best suited to
attack the isometric deformation problem via geometric transformation, since it provides correspondence of the
characteristics of the isometric deformation problem (according to Darboux these are the asymptotic
directions), it is directly linked to the infinitesimal isometric deformation problem (Darboux
proved that infinitesimal isometric deformations generate Weingarten congruences and Guichard proved
the converse: there is an infinitesimal isometric deformation of a focal surface of a Weingarten
congruence in the direction normal to the other focal surface; see Darboux (\cite{D1},\S\ 883-\S\
924)) and it admits a version of the Bianchi Permutability Theorem for its second iteration.

In \cite{D3} we proved that for a generic $3$-dimensional integrable rolling distribution of contact elements (excluding
developable seed and isotropic developable leaves) and with the symmetry of the tangency configuration (contact elements
are centered on tangent planes of the surface $x_0$ and further pass through the origin of the tangent planes) the seed
and any leaf are the focal surfaces of a Weingarten congruence (and thus we get B\"{a}cklund transformation according to
Bianchi's definition).

We have now the main {\bf Theorem} of this paper:

\begin{theorem}\label{th:th1}
A generic $4$-dimensional integrable rolling distribution of contact elements with the symmetry of the tangency
configuration (excluding developable seed and isotropic developable leaves) splits into an $1$-dimensional family
of generic $3$-dimensional integrable rolling distributions of contact elements with the symmetry of the tangency
configuration.
\end{theorem}

The remaining part of the paper is organized as follows: in Section 2 we recall the rolling problem for surfaces and in
Section 3 we provide the proof of {\bf Theorem} \ref{th:th1}.

\section{The rolling problem for surfaces}

Let $(u,v)\in D$ with $D$ domain of $\mathbb{R}^2$ or $\mathbb{C}^2$ and $x:D\mapsto\mathbb{C}^3$ be
a surface.

For $\om_1,\om_2\ \mathbb{C}^3$-valued $1$-forms on $D$ and $a,b\in\mathbb{C}^3$ we have
\begin{eqnarray}\label{eq:fund}
a^T\om_1\wedge b^T\om_2=((a\times b)\times\om_1+b^T\om_1a)^T\wedge\om_2=
(a\times b)^T(\om_1\times\wedge\om_2)+b^T\om_1\wedge a^T\om_2;\nonumber\\
\mathrm{in\ particular}\ a^T\om\wedge b^T\om=\frac{1}{2}(a\times b)^T(\om\times\wedge\om).
\end{eqnarray}

Since both $\times$ and $\wedge$ are skew-symmetric, we have $2\om_1\times\wedge\om_2=
\om_1\times\om_2+\om_2\times\om_1=2\om_2\times\wedge\om_1$.

Consider the scalar product $\langle\cdot,\cdot\rangle$ on $\mathbf{M}_3(\mathbb{C}):\
\langle X,Y\rangle :=\frac{1}{2}\mathrm{tr}(X^TY)$. We have the isometry

$$\al:\mathbb{C}^3\mapsto\mathbf{o}_3(\mathbb{C}),\
\al\left(\begin{bmatrix}x^1\\x^2\\x^3\end{bmatrix}\right)
=\begin{bmatrix}0&-x^3&x^2\\x^3&0&-x^1\\-x^2&x^1&0\end{bmatrix},\
x^Ty=\langle\al(x),\al(y)\rangle=
\frac{1}{2}\mathrm{tr}(\al(x)^T\al(y)),$$
$$\al(x\times y)=[\al(x),\al(y)]=\al(\al(x)y)=yx^T-xy^T,\ \al(Rx)=R\al(x)R^{-1},\
x,y\in\mathbb{C}^3,\ R\in\mathbf{O}_3(\mathbb{C}).$$

Let $x\subset\mathbb{C}^3$ be a surface applicable (isometric) to a surface
$x_0\subset\mathbb{C}^3$:

\begin{eqnarray}\label{eq:roll}
(x,dx)=(R,t)(x_0,dx_0):=(Rx_0+t,Rdx_0),
\end{eqnarray}
where $(R,t)$ is a sub-manifold in $\mathbf{O}_3(\mathbb{C})\ltimes\mathbb{C}^3$ (in general
surface, but it is a curve if $x_0,\ x$ are ruled and the rulings correspond under isometry or a
point if $x_0,\ x$ differ by a rigid motion). The sub-manifold $R$ gives the rolling of $x_0$ on $x$,
that is if we rigidly roll $x_0$ on $x$ such that points corresponding under the isometry will have
the same differentials, $R$ will dictate the rotation of $x_0$; the translation $t$ will satisfy
$dt=-dRx_0$.

For $(u,v)$ parametrization on $x_0,\ x$ and outside the locus of isotropic (degenerate) induced
metric of $x_0,\ x$ we have $N_0:=\frac{\pa_ux_0\times\pa_vx_0}{|\pa_ux_0\times\pa_vx_0|},\
N:=\frac{\pa_ux\times\pa_vx}{|\pa_ux\times\pa_vx|}$ respectively positively oriented unit normal
fields of $x_0,\ x$ and $R$ is determined by $R=[\pa_ux\ \ \pa_vx\ \ N][\pa_ux_0\ \ \pa_vx_0\ \
\det(R)N_0]^{-1}$; we take $R$ with $\det(R)=1$; thus the rotation of the rolling with the other
face of $x_0$ (or on the other face of $x$) is $R':=R(I_3-2N_0N_0^T)=(I_3-2NN^T)R,\ \det(R')=-1$.

Therefore $\mathbf{O}_3(\mathbb{C})\ltimes\mathbb{C}^3$ acts on
$2$-dimensional integrable distributions of contact elements $(x_0,dx_0)$ in
$T^*(\mathbb{C}^3)$ as: $(R,t)(x_0,dx_0)=(Rx_0+t,Rdx_0)$; a
rolling is a sub-manifold
$(R,t)\subset\mathbf{O}_3(\mathbb{C})\ltimes\mathbb{C}^3$ such
that $(R,t)(x_0,dx_0)$ is still integrable.

We have:
\begin{eqnarray}\label{eq:secoi}
R^{-1}dRN_0=R^{-1}dN-dN_0.
\end{eqnarray}
Applying the
compatibility condition $d$ to (\ref{eq:roll}) we get:
\begin{eqnarray}\label{eq:comp}
R^{-1}dR\wedge dx_0=0,\ dRR^{-1}\wedge dx=0.
\end{eqnarray}

Since $R^{-1}dR$ is skew-symmetric and using (\ref{eq:comp}) we
have
\begin{eqnarray}\label{eq:dx0}
dx_0^TR^{-1}dRdx_0=0.
\end{eqnarray}
From (\ref{eq:dx0}) for $a\in\mathbb{C}^3$ we get
$R^{-1}dRa=R^{-1}dR(a^{\bot}+a^{\top})=a^TN_0R^{-1}dRN_0-a^TR^{-1}dRN_0N_0=\om\times
a,\ \om:=N_0\times R^{-1}dRN_0=^{(\ref{eq:secoi})}(\det
R)R^{-1}(N\times dN)-N_0\times dN_0=R^{-1}(N\times dN)-N_0\times
dN_0$. Thus $R^{-1}dR=\al(\om)$ and $\om$ is flat connection form
in $T^*x_0$:
\begin{eqnarray}\label{eq:om}
d\om+\frac{1}{2}\om\times\wedge\om=0,\ \om\times\wedge
dx_0=0,\ (\om)^\perp=0.
\end{eqnarray}
With $s:=N_0^T(\om\times
dx_0)=s_{11}du^2+s_{12}dudv+s_{21}dvdu+s_{22}dv^2$ the difference
of the second fundamental forms of $x,\ x_0$ we have
\begin{eqnarray}\label{eq:omjk}
\om=\frac{s_{12}\pa_ux_0-s_{11}\pa_vx_0}{|\pa_ux_0\times\pa_vx_0|}du+
\frac{s_{22}\pa_ux_0-s_{21}\pa_vx_0}{|\pa_ux_0\times\pa_vx_0|}dv;
\end{eqnarray}
($\om\times\wedge dx_0=0$ is equivalent to $s_{12}=s_{21}$;
$(d\om)^\perp+\frac{1}{2}\om\times\wedge\om=0,\
(d\om)^\top=0$ respectively encode the difference of the
Gau\ss\ -Codazzi-Mainardi-Peterson equations of $x_0$ and $x$).

Using $\frac{1}{2}dN_0\times\wedge
dN_0=K|\pa_ux_0\times\pa_vx_0|N_0du\wedge dv,\ K$ being the
Gau\ss\ curvature we get $dN_0\times\wedge
dN_0=R^{-1}(dN\times\wedge dN)=^{(\ref{eq:secoi})}(\om\times
N_0+dN_0)\times\wedge(\om\times N_0+dN_0)=dN_0\times\wedge
dN_0+2(\om\times N_0)\times\wedge dN_0+\om\times\wedge\om$; thus
\begin{eqnarray}\label{eq:omom}
\frac{1}{2}\om\times\wedge\om=dN_0^T\wedge\om N_0.
\end{eqnarray}
Note also
\begin{eqnarray}\label{eq:om'}
\om'=N_0\times {R'}^{-1}dR'N_0=-\om-2N_0\times dN_0
\end{eqnarray}
and
\begin{eqnarray}\label{eq:aom}
\ \ \ \ a^T\wedge\om=0,\ \forall\om\ \mathrm{satisfying\
(\ref{eq:om})\ for}\ a\ 1-\mathrm{form}\Rightarrow\ a^T\odot
dx_0:=\frac{a^Tdx_0+dx_0^Ta}{2}=0.
\end{eqnarray}
Note that the converse $a^T\odot dx_0=0,\ a\ 1-$form
$\Rightarrow a^T\wedge\om=0,\forall\om$ satisfying (\ref{eq:om})
is also true.

\section {Proof of Theorem \ref{th:th1}}

Consider a surface $x_0=x_0(u,v)\subset\mathbb{C}^3$ with unit normal field $N_0=N_0(u,v)$.

Consider a $4$-dimensional distribution of contact elements with the symmetry of the tangency
configuration, that is the contact elements are centered at $x_0+V,\ V=V(u,v,w_1,w_2),\ N_0^TV=0,\
du\wedge dv\wedge dw_1\wedge dw_2\neq 0$ and have non-isotropic normal fields $m=V\times N_0+\mathbf{m}N_0,\
\mathbf{m}=\mathbf{m}(u,v,w_1,w_2)\subset\mathbb{C}$.

With $\ti d\cdot:=\pa_u\cdot du+\pa_v\cdot dv+\pa_{w_1}\cdot dw_1+\pa_{w_2}\cdot dw_2=d\cdot+\pa_{w_1}\cdot dw_1+
\pa_{w_2}\cdot dw_2$ if the
distribution of contact elements is integrable and the rolled distribution remains integrable if we
roll $x_0$ on an isometric surface $x,\ (x,dx)=(R,t)(x_0,dx_0)$ (that is we replace $x_0,\ V,\ m$
with $x,\ RV,\ Rm$), then along the leaves we have
$$0=(Rm)^T\ti d(RV+x)=m^T(\om\times V+d(V+x_0)+\pa_{w_1}Vdw_1+\pa_{w_2}Vdw_2),$$
or, assuming $N_0^T(\pa_{w_2}V\times V)\neq 0$,

\begin{eqnarray}\label{eq:dw2}
dw_2=\frac{N_0^T[V\times d(V+x_0)]}{N_0^T(\pa_{w_2}V\times V)}+
\mathbf{m}\frac{V^T(\om\times N_0+dN_0)}{N_0^T(\pa_{w_2}V\times V)}-
\frac{N_0^T(\pa_{w_1}V\times V)}{N_0^T(\pa_{w_2}V\times V)}dw_1.
\end{eqnarray}

By applying the compatibility condition $\ti d$ to (\ref{eq:dw2}) and using the equation itself we get
the integrability condition

$$0=\ti d(\frac{N_0^T[V\times d(V+x_0)]}{N_0^T(\pa_{w_2}V\times V)}+
\mathbf{m}\frac{V^T(\om\times N_0+dN_0)}{N_0^T(\pa_{w_2}V\times V)}-
\frac{N_0^T(\pa_{w_1}V\times V)}{N_0^T(\pa_{w_2}V\times V)}dw_1)$$
$$=-\frac{\pa_{w_2}[N_0^T[V\times d(V+x_0)]+\mathbf{m}V^T(\om\times N_0+dN_0)-
N_0^T(\pa_{w_1}V\times V)dw_1]}{N_0^T(\pa_{w_2}V\times V)}\wedge(\frac{N_0^T[V\times d(V+x_0)]}
{N_0^T(\pa_{w_2}V\times V)}+$$
$$\mathbf{m}\frac{V^T(\om\times N_0+dN_0)}{N_0^T(\pa_{w_2}V\times V)}-
\frac{N_0^T(\pa_{w_1}V\times V)}{N_0^T(\pa_{w_2}V\times V)}dw_1)-
\pa_{w_1}\frac{N_0^T[V\times d(V+x_0)]+\mathbf{m}V^T(\om\times N_0+dN_0)}{N_0^T(\pa_{w_2}V\times V)}\wedge dw_1$$
$$+d\frac{N_0^T[V\times d(V+x_0)]}{N_0^T(\pa_{w_2}V\times V)}+
d\frac{\mathbf{m}V^T}{N_0^T(\pa_{w_2}V\times V)}\wedge(\om\times N_0+dN_0)-
d\frac{N_0^T(\pa_{w_1}V\times V)}{N_0^T(\pa_{w_2}V\times V)}\wedge dw_1$$
$$\stackrel{(\ref{eq:fund})}{=}-
[\frac{N_0^T[\pa_{w_2}V\times d(V+x_0)]\wedge N_0^T(V\times dx_0)}{N_0^T(\pa_{w_2}V\times V)}+
(\mathbf{m}^2+|V|^2)\frac{1}{2}N_0^T(dN_0\times\wedge dN_0)]\frac{1}{N_0^T(\pa_{w_2}V\times V)}$$
$$-\frac{N_0^T(\pa_{w_1}V\times \pa_{w_2}V)N_0^T(V\times dx_0)}{N_0^T(\pa_{w_2}V\times V)^2}\wedge dw_1
+[d\mathbf{m}+\pa_{w_2}\mathbf{m}\frac{N_0^T[V\times d(V+x_0)]}{N_0^T(\pa_{w_2}V\times V)}-$$
$$\mathbf{m}\frac{N_0^T(\pa_{w_2}V\times dV)}{N_0^T(\pa_{w_2}V\times V)}+
(\mathbf{m}\frac{N_0^T(\pa_{w_1}V\times\pa_{w_2}V)}{N_0^T(\pa_{w_2}V\times V)}-
\pa_{w_2}\mathbf{m}\frac{N_0^T(\pa_{w_1}V\times V)}{N_0^T(\pa_{w_2}V\times V)}+
\pa_{w_1}\mathbf{m})dw_1]\wedge\frac{V^T(\om\times N_0+dN_0)}{N_0^T(\pa_{w_2}V\times V)}$$

for all $\om$ satisfying (\ref{eq:om}).

By exchanging the role played by $w_1$ and $w_2$ and using (\ref{eq:dw2}) the new integrability condition is the old one
multiplied by
$\frac{N_0^T(\pa_{w_2}V\times V)}{N_0^T(\pa_{w_1}V\times V)}$, so it essentially remains the same.

For our problem we should get compatibility conditions independent of the shape of the seed surface $x$ and the solution
should depend on two constants.

According to these principles and from (\ref{eq:dw2}) we make the ansatz that along the leaves we have

\begin{eqnarray}\label{eq:dw1}
dw_1=A_1+B_1V^T(\om\times N_0+dN_0),
\end{eqnarray}

where $A_1$ is a scalar $1$-form not depending on $\om,\ dw_1,\ dw_2$
and $B_1=B_1(u,v,w_1,w_2)\subset\mathbb{C}$ is a function.

With

$$\mathcal{A}_2:=[\frac{N_0^T[\pa_{w_2}V\times d(V+x_0)]\wedge N_0^T(V\times dx_0)}{N_0^T(\pa_{w_2}V\times V)}+
(\mathbf{m}^2+|V|^2)\frac{1}{2}N_0^T(dN_0\times\wedge dN_0)]\frac{1}{N_0^T(\pa_{w_2}V\times V)},$$
$$\mathcal{B}_2:=d\mathbf{m}+\pa_{w_2}\mathbf{m}\frac{N_0^T[V\times d(V+x_0)]}{N_0^T(\pa_{w_2}V\times V)}-
\mathbf{m}\frac{N_0^T(\pa_{w_2}V\times dV)}{N_0^T(\pa_{w_2}V\times V)},$$
$$\mathcal{C}_2:=\mathbf{m}\frac{N_0^T(\pa_{w_1}V\times\pa_{w_2}V)}{N_0^T(\pa_{w_2}V\times V)}-
\pa_{w_2}\mathbf{m}\frac{N_0^T(\pa_{w_1}V\times V)}{N_0^T(\pa_{w_2}V\times V)}+
\pa_{w_1}\mathbf{m}$$

from the integrability condition of (\ref{eq:dw2}) and using (\ref{eq:dw1}) we have

\begin{eqnarray}\label{eq:a1}
A_1=-\frac{\mathcal{B}_2}{\mathcal{C}_2}+
\frac{N_0^T(\pa_{w_1}V\times \pa_{w_2}V)N_0^T(V\times dx_0)}{N_0^T(\pa_{w_2}V\times V)}\frac{B_1}{\mathcal{C}_2},\nonumber\\
-\mathcal{A}_2+\frac{N_0^T(\pa_{w_1}V\times \pa_{w_2}V)N_0^T(V\times dx_0)}{N_0^T(\pa_{w_2}V\times V)^2}\wedge
\frac{\mathcal{B}_2}{\mathcal{C}_2}=0.
\end{eqnarray}

If we exchange the role played by $w_1$ and $w_2$ in the first equation of (\ref{eq:a1}) and using (\ref{eq:dw2}) and (\ref{eq:dw1}), then the new equation remains the same as the old one.

The second equation of (\ref{eq:a1}) is a compatibility condition for our problem independent of the shape of the seed
surface $x$; if we exchange the role played by $w_1$ and $w_2$, then the new corresponding equation will be the old one
multiplied by $\frac{N_0^T(\pa_{w_2}V\times V)}{N_0^T(\pa_{w_1}V\times V)}$, so it essentially remains the same.

From the first equation of (\ref{eq:a1}) and using (\ref{eq:dw1}) we get

\begin{eqnarray}\label{eq:dw1a}
dw_1=-\frac{\mathcal{B}_2}{\mathcal{C}_2}+B_1[\frac{N_0^T(\pa_{w_1}V\times \pa_{w_2}V)N_0^T(V\times dx_0)}
{N_0^T(\pa_{w_2}V\times V)\mathcal{C}_2}+V^T(\om\times N_0+dN_0)].
\end{eqnarray}

By applying the compatibility condition $\ti d$ to (\ref{eq:dw1a}) and using the equation itself and (\ref{eq:dw2})
we get the integrability condition

$$0=\ti d[-\frac{\mathcal{B}_2}{\mathcal{C}_2}+B_1[\frac{N_0^T(\pa_{w_1}V\times \pa_{w_2}V)N_0^T(V\times dx_0)}
{N_0^T(\pa_{w_2}V\times V)\mathcal{C}_2}+V^T(\om\times N_0+dN_0)]]=$$
$$-\pa_{w_1}[-\frac{\mathcal{B}_2}{\mathcal{C}_2}+B_1[\frac{N_0^T(\pa_{w_1}V\times \pa_{w_2}V)N_0^T(V\times dx_0)}
{N_0^T(\pa_{w_2}V\times V)\mathcal{C}_2}+V^T(\om\times N_0+dN_0)]]$$
$$\wedge(-\frac{\mathcal{B}_2}{\mathcal{C}_2}+B_1[\frac{N_0^T(\pa_{w_1}V\times \pa_{w_2}V)N_0^T(V\times dx_0)}
{N_0^T(\pa_{w_2}V\times V)\mathcal{C}_2}+V^T(\om\times N_0+dN_0)])-\pa_{w_2}[-\frac{\mathcal{B}_2}{\mathcal{C}_2}+$$
$$B_1[\frac{N_0^T(\pa_{w_1}V\times \pa_{w_2}V)N_0^T(V\times dx_0)}
{N_0^T(\pa_{w_2}V\times V)\mathcal{C}_2}+V^T(\om\times N_0+dN_0)]]\wedge(\frac{N_0^T[V\times d(V+x_0)]}
{N_0^T(\pa_{w_2}V\times V)}+$$
$$\mathbf{m}\frac{V^T(\om\times N_0+dN_0)}{N_0^T(\pa_{w_2}V\times V)}-\frac{N_0^T(\pa_{w_1}V\times V)}{N_0^T(\pa_{w_2}V\times V)}[-\frac{\mathcal{B}_2}{\mathcal{C}_2}
+B_1[\frac{N_0^T(\pa_{w_1}V\times \pa_{w_2}V)N_0^T(V\times dx_0)}
{N_0^T(\pa_{w_2}V\times V)\mathcal{C}_2}+$$
$$V^T(\om\times N_0+dN_0)]])+d(-\frac{\mathcal{B}_2}{\mathcal{C}_2}+B_1\frac{N_0^T(\pa_{w_1}V\times \pa_{w_2}V)N_0^T(V\times dx_0)}
{N_0^T(\pa_{w_2}V\times V)\mathcal{C}_2})+d(B_1V^T)\wedge(\om\times N_0+dN_0)$$
$$\stackrel{(\ref{eq:fund})}{=}\frac{d\pa_{w_1}\mathbf{m}}{\mathcal{C}_2}0+
\frac{d\pa_{w_2}\mathbf{m}}{\mathcal{C}_2}0+\frac{\pa_{w_1}^2\mathbf{m}}{\mathcal{C}_2}0
+\frac{\pa_{w_2}^2\mathbf{m}}{\mathcal{C}_2}0+\frac{\pa_{w_1w_2}^2\mathbf{m}}{\mathcal{C}_2}0+$$
$$\frac{d\mathbf{m}}{\mathcal{C}_2}\wedge[\frac{\pa_{w_2}\mathbf{m}\pa_{w_1}\frac{N_0^T[V\times d(V+x_0)]}{N_0^T(\pa_{w_2}V\times V)}-
\pa_{w_1}(\mathbf{m}\frac{N_0^T(\pa_{w_2}V\times dV)}{N_0^T(\pa_{w_2}V\times V)})}{\mathcal{C}_2}-\frac{1}{\mathcal{C}_2}\pa_{w_1}\frac{B_1N_0^T(\pa_{w_1}V\times \pa_{w_2}V)N_0^T(V\times dx_0)}{N_0^T(\pa_{w_2}V\times V)}$$$$-\frac{\pa_{w_2}(\mathbf{m}\frac{N_0^T(\pa_{w_1}V\times\pa_{w_2}V)}{N_0^T(\pa_{w_2}V\times V)})-
\pa_{w_2}\mathbf{m}\pa_{w_2}\frac{N_0^T(\pa_{w_1}V\times V)}{N_0^T(\pa_{w_2}V\times V)}}{\mathcal{C}_2}
\frac{N_0^T[V\times d(V+x_0)]}{N_0^T(\pa_{w_2}V\times V)}-\frac{N_0^T(\pa_{w_1}V\times V)}{N_0^T(\pa_{w_2}V\times V)}$$$$(\frac{\pa_{w_2}\mathbf{m}\pa_{w_2}\frac{N_0^T[V\times d(V+x_0)]}{N_0^T(\pa_{w_2}V\times V)}-
\pa_{w_2}(\mathbf{m}\frac{N_0^T(\pa_{w_2}V\times dV)}{N_0^T(\pa_{w_2}V\times V)})}{\mathcal{C}_2}-\frac{1}{\mathcal{C}_2}\pa_{w_2}\frac{B_1N_0^T(\pa_{w_1}V\times \pa_{w_2}V)N_0^T(V\times dx_0)}{N_0^T(\pa_{w_2}V\times V)})-$$$$\frac{\mathbf{m}d\frac{N_0^T(\pa_{w_1}V\times\pa_{w_2}V)}{N_0^T(\pa_{w_2}V\times V)}-
\pa_{w_2}\mathbf{m}d\frac{N_0^T(\pa_{w_1}V\times V)}{N_0^T(\pa_{w_2}V\times V)}}{\mathcal{C}_2}+
\frac{N_0^T(\pa_{w_2}V\times dV)}{N_0^T(\pa_{w_2}V\times V)}-\frac{N_0^T(\pa_{w_1}V\times\pa_{w_2}V)}{N_0^T(\pa_{w_2}V\times V)}B_1$$
$$\frac{N_0^T(\pa_{w_1}V\times \pa_{w_2}V)N_0^T(V\times dx_0)}
{N_0^T(\pa_{w_2}V\times V)\mathcal{C}_2}]+\frac{\pa_{w_1}\mathbf{m}}{\mathcal{C}_2}[\frac{N_0^T(\pa_{w_2}V\times dV)}{N_0^T(\pa_{w_2}V\times V)}\wedge(\pa_{w_2}\mathbf{m}\frac{N_0^T[V\times d(V+x_0)]}{N_0^T(\pa_{w_2}V\times V)\mathcal{C}_2}-$$$$B_1\frac{N_0^T(\pa_{w_1}V\times \pa_{w_2}V)N_0^T(V\times dx_0)}
{N_0^T(\pa_{w_2}V\times V)\mathcal{C}_2})-\frac{\pa_{w_2}\mathbf{m}\frac{N_0^T[V\times d(V+x_0)]}{N_0^T(\pa_{w_2}V\times V)}-
\mathbf{m}\frac{N_0^T(\pa_{w_2}V\times dV)}{N_0^T(\pa_{w_2}V\times V)}}{\mathcal{C}_2}\frac{N_0^T(\pa_{w_1}V\times\pa_{w_2}V)}{N_0^T(\pa_{w_2}V\times V)}$$$$\wedge B_1\frac{N_0^T(\pa_{w_1}V\times \pa_{w_2}V)N_0^T(V\times dx_0)}
{N_0^T(\pa_{w_2}V\times V)\mathcal{C}_2}-B_1\frac{N_0^T(\pa_{w_1}V\times \pa_{w_2}V)N_0^T(V\times dx_0)}
{N_0^T(\pa_{w_2}V\times V)^2\mathcal{C}_2}N_0^T(\pa_{w_1w_2}^2V\times V)\wedge$$$$\frac{\pa_{w_2}\mathbf{m}\frac{N_0^T[V\times d(V+x_0)]}{N_0^T(\pa_{w_2}V\times V)}-
\mathbf{m}\frac{N_0^T(\pa_{w_2}V\times dV)}{N_0^T(\pa_{w_2}V\times V)}}{\mathcal{C}_2}+B_1\frac{N_0^T(\pa_{w_1}V\times \pa_{w_2}V)N_0^T(V\times dx_0)}
{N_0^T(\pa_{w_2}V\times V)^2\mathcal{C}_2}N_0^T(\pa_{w_2}^2V\times V)\wedge$$
$$(\frac{N_0^T[V\times d(V+x_0)]}{N_0^T(\pa_{w_2}V\times V)}+\frac{N_0^T(\pa_{w_1}V\times V)}{N_0^T(\pa_{w_2}V\times V)}\frac{\pa_{w_2}\mathbf{m}\frac{N_0^T[V\times d(V+x_0)]}{N_0^T(\pa_{w_2}V\times V)}-
\mathbf{m}\frac{N_0^T(\pa_{w_2}V\times dV)}{N_0^T(\pa_{w_2}V\times V)}}{\mathcal{C}_2})-$$$$\frac{d(N_0^T(\pa_{w_2}V\times V))}{N_0^T(\pa_{w_2}V\times V)}\wedge B_1\frac{N_0^T(\pa_{w_1}V\times \pa_{w_2}V)N_0^T(V\times dx_0)}
{N_0^T(\pa_{w_2}V\times V)\mathcal{C}_2}]+\frac{\pa_{w_2}\mathbf{m}}{\mathcal{C}_2}[\pa_{w_1}\frac{N_0^T[V\times d(V+x_0)]}{N_0^T(\pa_{w_2}V\times V)}\wedge$$$$(-\frac{\pa_{w_2}\mathbf{m}\frac{N_0^T[V\times d(V+x_0)]}{N_0^T(\pa_{w_2}V\times V)}-
\mathbf{m}\frac{N_0^T(\pa_{w_2}V\times dV)}{N_0^T(\pa_{w_2}V\times V)}}{\mathcal{C}_2}+B_1\frac{N_0^T(\pa_{w_1}V\times \pa_{w_2}V)N_0^T(V\times dx_0)}
{N_0^T(\pa_{w_2}V\times V)\mathcal{C}_2})+$$$$\frac{\pa_{w_2}\mathbf{m}\frac{N_0^T[V\times d(V+x_0)]}{N_0^T(\pa_{w_2}V\times V)}-
\mathbf{m}\frac{N_0^T(\pa_{w_2}V\times dV)}{N_0^T(\pa_{w_2}V\times V)}}{\mathcal{C}_2}\pa_{w_1}\frac{N_0^T(\pa_{w_1}V\times V)}{N_0^T(\pa_{w_2}V\times V)}\wedge B_1\frac{N_0^T(\pa_{w_1}V\times \pa_{w_2}V)N_0^T(V\times dx_0)}
{N_0^T(\pa_{w_2}V\times V)\mathcal{C}_2}$$$$+B_1\frac{N_0^T(\pa_{w_1}V\times \pa_{w_2}V)N_0^T(V\times dx_0)}
{N_0^T(\pa_{w_2}V\times V)^2\mathcal{C}_2}\pa_{w_1}(N_0^T(\pa_{w_1}V\times V))\wedge\frac{\pa_{w_2}\mathbf{m}\frac{N_0^T[V\times d(V+x_0)]}{N_0^T(\pa_{w_2}V\times V)}-
\mathbf{m}\frac{N_0^T(\pa_{w_2}V\times dV)}{N_0^T(\pa_{w_2}V\times V)}}{\mathcal{C}_2}$$$$+\frac{\pa_{w_1}[B_1N_0^T(\pa_{w_1}V\times \pa_{w_2}V)N_0^T(V\times dx_0)]}
{N_0^T(\pa_{w_2}V\times V)\mathcal{C}_2}\wedge\frac{N_0^T[V\times d(V+x_0)]}{N_0^T(\pa_{w_2}V\times V)}$$
$$+\frac{\mathbf{m}}{\mathcal{C}_2}\pa_{w_1}\frac{N_0^T(\pa_{w_2}V\times dV)}{N_0^T(\pa_{w_2}V\times V)}\wedge\frac{N_0^T[V\times d(V+x_0)]}{N_0^T(\pa_{w_2}V\times V)}+(\pa_{w_2}\frac{N_0^T[V\times d(V+x_0)]}{N_0^T(\pa_{w_2}V\times V)}$$$$-\frac{N_0^T(\pa_{w_2}V\times dV)}{N_0^T(\pa_{w_2}V\times V)})\wedge[\frac{N_0^T[V\times d(V+x_0)]}{N_0^T(\pa_{w_2}V\times V)}-\frac{N_0^T(\pa_{w_1}V\times V)}{N_0^T(\pa_{w_2}V\times V)}(-\frac{\pa_{w_2}\mathbf{m}\frac{N_0^T[V\times d(V+x_0)]}{N_0^T(\pa_{w_2}V\times V)}-
\mathbf{m}\frac{N_0^T(\pa_{w_2}V\times dV)}{N_0^T(\pa_{w_2}V\times V)}}{\mathcal{C}_2}$$$$+B_1\frac{N_0^T(\pa_{w_1}V\times \pa_{w_2}V)N_0^T(V\times dx_0)}
{N_0^T(\pa_{w_2}V\times V)\mathcal{C}_2})]+(-\frac{N_0^T(\pa_{w_1}V\times\pa_{w_2}V)}{N_0^T(\pa_{w_2}V\times V)}+\pa_{w_2}\frac{N_0^T(\pa_{w_1}V\times V)}{N_0^T(\pa_{w_2}V\times V)})$$$$\frac{\pa_{w_2}\mathbf{m}\frac{N_0^T[V\times d(V+x_0)]}{N_0^T(\pa_{w_2}V\times V)}-
\mathbf{m}\frac{N_0^T(\pa_{w_2}V\times dV)}{N_0^T(\pa_{w_2}V\times V)}}{\mathcal{C}_2}\wedge[\frac{N_0^T[V\times d(V+x_0)]}{N_0^T(\pa_{w_2}V\times V)}-\frac{N_0^T(\pa_{w_1}V\times V)}{N_0^T(\pa_{w_2}V\times V)}B_1$$$$\frac{N_0^T(\pa_{w_1}V\times \pa_{w_2}V)N_0^T(V\times dx_0)}
{N_0^T(\pa_{w_2}V\times V)\mathcal{C}_2}]-B_1\frac{N_0^T(\pa_{w_1}V\times \pa_{w_2}V)N_0^T(V\times dx_0)}
{N_0^T(\pa_{w_2}V\times V)^2\mathcal{C}_2}N_0^T(\pa_{w_1w_2}^2V\times V)$$$$\wedge[\frac{N_0^T[V\times d(V+x_0)]}{N_0^T(\pa_{w_2}V\times V)}+\frac{N_0^T(\pa_{w_1}V\times V)}{N_0^T(\pa_{w_2}V\times V)}\frac{\pa_{w_2}\mathbf{m}\frac{N_0^T[V\times d(V+x_0)]}{N_0^T(\pa_{w_2}V\times V)}-\mathbf{m}\frac{N_0^T(\pa_{w_2}V\times dV)}{N_0^T(\pa_{w_2}V\times V)}}{\mathcal{C}_2}]$$$$-\frac{\pa_{w_2}[B_1N_0^T(\pa_{w_1}V\times \pa_{w_2}V)N_0^T(V\times dx_0)]}
{N_0^T(\pa_{w_2}V\times V)\mathcal{C}_2}\wedge\frac{N_0^T(\pa_{w_1}V\times V)}{N_0^T(\pa_{w_2}V\times V)}\frac{N_0^T[V\times d(V+x_0)]}{N_0^T(\pa_{w_2}V\times V)}$$$$-\frac{\mathbf{m}}{\mathcal{C}_2}\pa_{w_2}\frac{N_0^T(\pa_{w_2}V\times dV)}{N_0^T(\pa_{w_2}V\times V)}\wedge\frac{N_0^T[V\times d(V+x_0)]}{N_0^T(\pa_{w_2}V\times V)}$$$$-d\frac{N_0^T[V\times d(V+x_0)]}{N_0^T(\pa_{w_2}V\times V)}+\frac{\pa_{w_2}\mathbf{m}\frac{N_0^T[V\times d(V+x_0)]}{N_0^T(\pa_{w_2}V\times V)}-
\mathbf{m}\frac{N_0^T(\pa_{w_2}V\times dV)}{N_0^T(\pa_{w_2}V\times V)}}{\mathcal{C}_2}\wedge d\frac{N_0^T(\pa_{w_1}V\times V)}{N_0^T(\pa_{w_2}V\times V)}-$$$$B_1\frac{N_0^T(\pa_{w_1}V\times \pa_{w_2}V)N_0^T(V\times dx_0)}
{N_0^T(\pa_{w_2}V\times V)^2\mathcal{C}_2}\wedge d(N_0^T(\pa_{w_1}V\times V))]$$$$+\frac{\mathbf{m}}{\mathcal{C}_2}[(-\pa_{w_1}\frac{N_0^T(\pa_{w_2}V\times dV)}{N_0^T(\pa_{w_2}V\times V)}+\mathbf{m}\frac{N_0^T(\pa_{w_2}V\times dV)}{\mathcal{C}_2N_0^T(\pa_{w_2}V\times V)}\pa_{w_1}\frac{N_0^T(\pa_{w_1}V\times\pa_{w_2}V)}{N_0^T(\pa_{w_2}V\times V)}+$$$$B_1\frac{N_0^T(\pa_{w_1}V\times \pa_{w_2}V)N_0^T(V\times dx_0)}
{N_0^T(\pa_{w_2}V\times V)^2\mathcal{C}_2}\pa_{w_1}N_0^T(\pa_{w_1}V\times \pa_{w_2}V))\wedge(\mathbf{m}\frac{N_0^T(\pa_{w_2}V\times dV)}{\mathcal{C}_2N_0^T(\pa_{w_2}V\times V)}+$$$$B_1\frac{N_0^T(\pa_{w_1}V\times \pa_{w_2}V)N_0^T(V\times dx_0)}
{N_0^T(\pa_{w_2}V\times V)\mathcal{C}_2})-\frac{\pa_{w_1}[B_1N_0^T(\pa_{w_1}V\times\pa_{w_2}V)N_0^T(V\times dx_0)]}{N_0^T(\pa_{w_2}V\times V)}\wedge\frac{N_0^T(\pa_{w_2}V\times dV)}{\mathcal{C}_2N_0^T(\pa_{w_2}V\times V)}$$$$+(-\pa_{w_2}\frac{N_0^T(\pa_{w_2}V\times dV)}{N_0^T(\pa_{w_2}V\times V)}+\mathbf{m}\frac{N_0^T(\pa_{w_2}V\times dV)}{\mathcal{C}_2N_0^T(\pa_{w_2}V\times V)}\pa_{w_2}\frac{N_0^T(\pa_{w_1}V\times\pa_{w_2}V)}{N_0^T(\pa_{w_2}V\times V)}+$$$$B_1\frac{N_0^T(\pa_{w_1}V\times \pa_{w_2}V)N_0^T(V\times dx_0)}
{N_0^T(\pa_{w_2}V\times V)^2\mathcal{C}_2}\pa_{w_2}N_0^T(\pa_{w_1}V\times \pa_{w_2}V))\wedge(\frac{N_0^T[V\times d(V+x_0)]}{N_0^T(\pa_{w_2}V\times V)}-\frac{N_0^T(\pa_{w_1}V\times V)}{N_0^T(\pa_{w_2}V\times V)}$$$$(\frac{\mathbf{m}N_0^T(\pa_{w_2}V\times dV)}{N_0^T(\pa_{w_2}V\times V)\mathcal{C}_2}+B_1\frac{N_0^T(\pa_{w_1}V\times \pa_{w_2}V)N_0^T(V\times dx_0)}
{N_0^T(\pa_{w_2}V\times V)\mathcal{C}_2}))+\frac{\pa_{w_2}[B_1N_0^T(\pa_{w_1}V\times \pa_{w_2}V)N_0^T(V\times dx_0)]}
{N_0^T(\pa_{w_2}V\times V)}$$$$\wedge\frac{N_0^T(\pa_{w_1}V\times V)}{N_0^T(\pa_{w_2}V\times V)}\frac{N_0^T(\pa_{w_2}V\times dV)}{N_0^T(\pa_{w_2}V\times V)\mathcal{C}_2}+d\frac{N_0^T(\pa_{w_2}V\times dV)}{N_0^T(\pa_{w_2}V\times V)}+\mathbf{m}\frac{N_0^T(\pa_{w_2}V\times dV)}{N_0^T(\pa_{w_2}V\times V)}\wedge\frac{1}{\mathcal{C}_2}d\frac{N_0^T(\pa_{w_1}V\times\pa_{w_2}V)}{N_0^T(\pa_{w_2}V\times V)}$$$$+B_1\frac{N_0^T(\pa_{w_1}V\times \pa_{w_2}V)N_0^T(V\times dx_0)}
{N_0^T(\pa_{w_2}V\times V)^2\mathcal{C}_2}\wedge d(N_0^T(\pa_{w_1}V\times \pa_{w_2}V))]-$$$$
\frac{\pa_{w_2}[B_1N_0^T(\pa_{w_1}V\times \pa_{w_2}V)N_0^T(V\times dx_0)]}
{N_0^T(\pa_{w_2}V\times V)\mathcal{C}_2}\wedge\frac{N_0^T[V\times d(V+x_0)]}{N_0^T(\pa_{w_2}V\times V)}-\frac{1}{2}\mathbf{m}B_1N_0^T(dN_0\times\wedge dN_0)+$$$$\frac{d[B_1N_0^T(\pa_{w_1}V\times \pa_{w_2}V)N_0^T(V\times dx_0)]}
{N_0^T(\pa_{w_2}V\times V)\mathcal{C}_2}$$

$$+[\frac{d\pa_{w_1}\mathbf{m}}{\mathcal{C}_2}B_1+\frac{d\pa_{w_2}\mathbf{m}}{\mathcal{C}_2}
(\frac{\mathbf{m}}{N_0^T(\pa_{w_2}V\times V)}-B_1\frac{N_0^T(\pa_{w_1}V\times V)}{N_0^T(\pa_{w_2}V\times V)})-\frac{\pa_{w_1}^2\mathbf{m}}{\mathcal{C}_2}\frac{\mathcal{B}_2}{\mathcal{C}_2}B_1+$$
$$\frac{\pa_{w_2}^2\mathbf{m}}{\mathcal{C}_2}[\frac{N_0^T[V\times d(V+x_0)]}{N_0^T(\pa_{w_2}V\times V)}+\frac{N_0^T(\pa_{w_1}V\times V)}{N_0^T(\pa_{w_2}V\times V)}(\frac{\mathcal{B}_2}{\mathcal{C}_2}-B_1\frac{N_0^T(\pa_{w_1}V\times \pa_{w_2}V)N_0^T(V\times dx_0)}
{N_0^T(\pa_{w_2}V\times V)\mathcal{C}_2})]$$
$$(\frac{\mathbf{m}}{N_0^T(\pa_{w_2}V\times V)}-B_1\frac{N_0^T(\pa_{w_1}V\times V)}{N_0^T(\pa_{w_2}V\times V)})+\frac{\pa_{w_1w_2}^2\mathbf{m}}{\mathcal{C}_2}[(\frac{N_0^T[V\times d(V+x_0)]}{N_0^T(\pa_{w_2}V\times V)}+\frac{N_0^T(\pa_{w_1}V\times V)}{N_0^T(\pa_{w_2}V\times V)}(\frac{\mathcal{B}_2}{\mathcal{C}_2}-$$$$B_1\frac{N_0^T(\pa_{w_1}V\times \pa_{w_2}V)N_0^T(V\times dx_0)}
{N_0^T(\pa_{w_2}V\times V)\mathcal{C}_2}))B_1-(\frac{\mathcal{B}_2}{\mathcal{C}_2}-B_1\frac{N_0^T(\pa_{w_1}V\times \pa_{w_2}V)N_0^T(V\times dx_0)}
{N_0^T(\pa_{w_2}V\times V)\mathcal{C}_2})$$$$(\frac{\mathbf{m}}{N_0^T(\pa_{w_2}V\times V)}-B_1\frac{N_0^T(\pa_{w_1}V\times V)}{N_0^T(\pa_{w_2}V\times V)})]
+\frac{d\mathbf{m}}{\mathcal{C}_2}[-\frac{\pa_{w_1}(\mathbf{m}\frac{N_0^T(\pa_{w_1}V
\times\pa_{w_2}V)}{N_0^T(\pa_{w_2}V\times V)})-
\pa_{w_2}\mathbf{m}\pa_{w_1}\frac{N_0^T(\pa_{w_1}V\times V)}{N_0^T(\pa_{w_2}V\times V)}}{\mathcal{C}_2}B_1$$$$-\pa_{w_1}B_1-\frac{\pa_{w_2}(\mathbf{m}\frac{N_0^T(\pa_{w_1}V
\times\pa_{w_2}V)}{N_0^T(\pa_{w_2}V\times V)})-
\pa_{w_2}\mathbf{m}\pa_{w_2}\frac{N_0^T(\pa_{w_1}V\times V)}{N_0^T(\pa_{w_2}V\times V)}}{\mathcal{C}_2}
(\frac{\mathbf{m}}{N_0^T(\pa_{w_2}V\times V)}-B_1\frac{N_0^T(\pa_{w_1}V\times V)}{N_0^T(\pa_{w_2}V\times V)})$$$$+\pa_{w_2}B_1\frac{N_0^T(\pa_{w_1}V\times V)}{N_0^T(\pa_{w_2}V\times V)}]
+\frac{\pa_{w_1}\mathbf{m}}{\mathcal{C}_2}[-B_1\frac{N_0^T(\pa_{w_2}V\times dV)}{N_0^T(\pa_{w_2}V\times V)}-B_1\frac{N_0^T(\pa_{w_1}V\times\pa_{w_2}V)}{N_0^T(\pa_{w_2}V\times V)}$$$$\frac{\pa_{w_2}\mathbf{m}\frac{N_0^T[V\times d(V+x_0)]}{N_0^T(\pa_{w_2}V\times V)}-
\mathbf{m}\frac{N_0^T(\pa_{w_2}V\times dV)}{N_0^T(\pa_{w_2}V\times V)}}{\mathcal{C}_2}+B_1^2\frac{N_0^T(\pa_{w_1}V\times \pa_{w_2}V)N_0^T(V\times dx_0)}
{N_0^T(\pa_{w_2}V\times V)^2\mathcal{C}_2}N_0^T(\pa_{w_1w_2}^2V\times V)$$$$+B_1\frac{N_0^T(\pa_{w_1}V\times \pa_{w_2}V)N_0^T(V\times dx_0)}
{N_0^T(\pa_{w_2}V\times V)^2\mathcal{C}_2}N_0^T(\pa_{w_2}^2V\times V)(\frac{\mathbf{m}}{N_0^T(\pa_{w_2}V\times V)}-B_1\frac{N_0^T(\pa_
{w_1}V\times V)}{N_0^T(\pa_{w_2}V\times V)})]+$$$$\frac{\pa_{w_2}\mathbf{m}}{\mathcal{C}_2}[B_1(\pa_{w_1}\frac{N_0^T[V\times d(V+x_0)]}{N_0^T(\pa_{w_2}V\times V)}+\frac{\pa_{w_2}\mathbf{m}\frac{N_0^T[V\times d(V+x_0)]}{N_0^T(\pa_{w_2}V\times V)}-
\mathbf{m}\frac{N_0^T(\pa_{w_2}V\times dV)}{N_0^T(\pa_{w_2}V\times V)}}{\mathcal{C}_2}\pa_{w_1}\frac{N_0^T(\pa_
{w_1}V\times V)}{N_0^T(\pa_{w_2}V\times V)}-$$$$B_1\frac{N_0^T(\pa_{w_1}V\times \pa_{w_2}V)N_0^T(V\times dx_0)}
{N_0^T(\pa_{w_2}V\times V)^2\mathcal{C}_2}N_0^T(\pa_{w_1}^2V\times V))-\pa_{w_1}B_1\frac{N_0^T[V\times d(V+x_0)]}{N_0^T(\pa_{w_2}V\times V)}+$$$$(\pa_{w_2}\frac{N_0^T[V\times d(V+x_0)]}{N_0^T(\pa_{w_2}V\times V)}-\frac{N_0^T(\pa_{w_2}V\times dV)}{N_0^T(\pa_{w_2}V\times V)}-(\frac{\pa_{w_2}\mathbf{m}\frac{N_0^T[V\times d(V+x_0)]}{N_0^T(\pa_{w_2}V\times V)}-\mathbf{m}\frac{N_0^T(\pa_{w_2}V\times dV)}{N_0^T(\pa_{w_2}V\times V)}}{\mathcal{C}_2}-$$$$B_1\frac{N_0^T(\pa_{w_1}V\times \pa_{w_2}V)N_0^T(V\times dx_0)}
{N_0^T(\pa_{w_2}V\times V)\mathcal{C}_2})(\frac{N_0^T(\pa_{w_1}V\times\pa_{w_2}V)}{N_0^T(\pa_{w_2}V\times V)}-\pa_{w_2}\frac{N_0^T(\pa_{w_1}V\times V)}{N_0^T(\pa_{w_2}V\times V)}))(\frac{\mathbf{m}}{N_0^T(\pa_{w_2}V\times V)}$$$$-B_1\frac{N_0^T(\pa_{w_1}V\times V)}{N_0^T(\pa_{w_2}V\times V)})+\pa_{w_2}B_1\frac{N_0^T[V\times d(V+x_0)]}{N_0^T(\pa_{w_2}V\times V)}\frac{N_0^T(\pa_{w_1}V\times V)}{N_0^T(\pa_{w_2}V\times V)}]+\frac{\mathbf{m}}{\mathcal{C}_2}[(-\pa_{w_1}\frac{N_0^T(\pa_{w_2}V\times dV)}{N_0^T(\pa_{w_2}V\times V)}+$$$$\mathbf{m}\frac{N_0^T(\pa_{w_2}V\times dV)}{\mathcal{C}_2N_0^T(\pa_{w_2}V\times V)}\pa_{w_1}\frac{N_0^T(\pa_{w_1}V\times\pa_{w_2}V)}{N_0^T(\pa_{w_2}V\times V)}+B_1\frac{N_0^T(\pa_{w_1}V\times \pa_{w_2}V)N_0^T(V\times dx_0)}
{N_0^T(\pa_{w_2}V\times V)^2\mathcal{C}_2}$$$$\pa_{w_1}N_0^T(\pa_{w_1}V\times \pa_{w_2}V)+\pa_{w_1}B_1\frac{N_0^T(\pa_{w_2}V\times dV)}{N_0^T(\pa_{w_2}V\times V)})B_1+(-\pa_{w_2}\frac{N_0^T(\pa_{w_2}V\times dV)}{N_0^T(\pa_{w_2}V\times V)}+\mathbf{m}\frac{N_0^T(\pa_{w_2}V\times dV)}{\mathcal{C}_2N_0^T(\pa_{w_2}V\times V)}$$$$\pa_{w_2}\frac{N_0^T(\pa_{w_1}V\times\pa_{w_2}V)}{N_0^T(\pa_{w_2}V\times V)}+B_1\frac{N_0^T(\pa_{w_1}V\times \pa_{w_2}V)N_0^T(V\times dx_0)}
{N_0^T(\pa_{w_2}V\times V)^2\mathcal{C}_2}\pa_{w_2}N_0^T(\pa_{w_1}V\times \pa_{w_2}V))$$$$(\frac{\mathbf{m}}{N_0^T(\pa_{w_2}V\times V)}-B_1\frac{N_0^T(\pa_{w_1}V\times V)}{N_0^T(\pa_{w_2}V\times V)})-\frac{\pa_{w_2}[B_1N_0^T(\pa_{w_1}V\times \pa_{w_2}V)N_0^T(V\times dx_0)]}
{N_0^T(\pa_{w_2}V\times V)^2}$$$$-\pa_{w_2}B_1\frac{N_0^T(\pa_{w_1}V\times V)}{N_0^T(\pa_{w_2}V\times V)}\frac{N_0^T(\pa_{w_2}V\times dV)}{N_0^T(\pa_{w_2}V\times V)}]-B_1^2\frac{\pa_{w_1}[N_0^T(\pa_{w_1}V\times \pa_{w_2}V)N_0^T(V\times dx_0)]}
{N_0^T(\pa_{w_2}V\times V)\mathcal{C}_2}+$$$$\pa_{w_2}B_1(\frac{N_0^T[V\times d(V+x_0)]}{N_0^T(\pa_{w_2}V\times V)}-\frac{N_0^T(\pa_{w_1}V\times \pa_{w_2}V)N_0^T(V\times dx_0)}{N_0^T(\pa_{w_2}V\times V)\mathcal{C}_2}\frac{\mathbf{m}}{N_0^T(\pa_{w_2}V\times V)})-$$$$B_1\frac{\pa_{w_2}[N_0^T(\pa_{w_1}V\times \pa_{w_2}V)N_0^T(V\times dx_0)]}{N_0^T(\pa_{w_2}V\times V)\mathcal{C}_2}(\frac{\mathbf{m}}{N_0^T(\pa_{w_2}V\times V)}-B_1\frac{N_0^T(\pa_{w_1}V\times V)}{N_0^T(\pa_{w_2}V\times V)})+dB_1]$$$$\wedge V^T(\om\times N_0+dN_0)
+B_1(-\frac{\mathcal{B}_2}{\mathcal{C}_2}+B_1\frac{N_0^T(\pa_{w_1}V\times \pa_{w_2}V)N_0^T(V\times dx_0)}{N_0^T(\pa_{w_2}V\times V)\mathcal{C}_2})\wedge\pa_{w_1}V^T(\om\times N_0+dN_0)$$
$$+B_1[\frac{N_0^T[V\times d(V+x_0)]}{N_0^T(\pa_{w_2}V\times V)}-\frac{N_0^T(\pa_{w_1}V\times V)}{N_0^T(\pa_{w_2}V\times V)}(-\frac{\mathcal{B}_2}{\mathcal{C}_2}+B_1\frac{N_0^T(\pa_{w_1}V\times \pa_{w_2}V)N_0^T(V\times dx_0)}{N_0^T(\pa_{w_2}V\times V)\mathcal{C}_2})]\wedge$$$$\pa_{w_2}V^T(\om\times N_0+dN_0)+B_1dV^T\wedge(\om\times N_0+dN_0)$$

for all $\om$ satisfying (\ref{eq:om}).

The terms containing $\pa_{w_1}V^T(\om\times N_0+dN_0),\pa_{w_2}V^T(\om\times N_0+dN_0)$ and $dV^T\wedge(\om\times N_0+dN_0)$ become $$[-B_1(-\frac{\mathcal{B}_2}{\mathcal{C}_2}+B_1\frac{N_0^T(\pa_{w_1}V\times \pa_{w_2}V)N_0^T(V\times dx_0)}{N_0^T(\pa_{w_2}V\times V)\mathcal{C}_2})\frac{N_0^T(\pa_{w_1}V\times \pa_{w_2}V)}{N_0^T(\pa_{w_2}V\times V)}+B_1\frac{N_0^T[\pa_{w_2}V\times d(V+x_0)]}{N_0^T(\pa_{w_2}V\times V)}]$$$$\wedge V^T(\om\times N_0+dN_0)$$
(use $dx_0^T\wedge(\om\times N_0+dN_0)=0$).

After further simplifications and using (\ref{eq:fund}) the integrability condition becomes

$$0=\frac{d\mathbf{m}}{\mathcal{C}_2}\wedge[-\frac{\mathbf{m}}{\mathcal{C}_2}(\frac{N_0^T(\pa_{w_1w_2}^2V\times \pa_{w_2}V)}{N_0^T(\pa_{w_2}V\times V)}-\frac{N_0^T(\pa_{w_1}V\times V)}{N_0^T(\pa_{w_2}V\times V)}\frac{N_0^T(\pa_{w_2}^2V\times \pa_{w_2}V)}{N_0^T(\pa_{w_2}V\times V)})\frac{N_0^T(V\times dx_0)}{N_0^T(\pa_{w_2}V\times V)}+$$$$\frac{\mathbf{m}}{\mathcal{C}_2}\frac{N_0^T(\pa_{w_1}V\times\pa_{w_2}V)}{N_0^T(\pa_{w_2}V\times V)}\frac{N_0^T(\pa_{w_2}V\times dV)}{N_0^T(\pa_{w_2}V\times V)}-\frac{\pa_{w_2}\mathbf{m}}{\mathcal{C}_2}
\frac{N_0^T(\pa_{w_1}V\times\pa_{w_2}V)}{N_0^T(\pa_{w_2}V\times V)}\frac{N_0^T(V\times dV)}{N_0^T(\pa_{w_2}V\times V)}-\frac{N_0^T(V\times dx_0)}{\mathcal{C}_2}$$$$(\pa_{w_1}B_1\frac{N_0^T(\pa_{w_1}V\times \pa_{w_2}V)}{N_0^T(\pa_{w_2}V\times V)}+B_1\frac{N_0^T(\pa_{w_1}^2V\times\pa_{w_2}V)}{N_0^T(\pa_{w_2}V\times V)}+B_1\frac{N_0^T(\pa_{w_1}V\times \pa_{w_2}V)^2}{N_0^T(\pa_{w_2}V\times V)^2}-\pa_{w_2}B_1\frac{N_0^T(\pa_{w_1}V\times V)}{N_0^T(\pa_{w_2}V\times V)}$$$$\frac{N_0^T(\pa_{w_1}V\times \pa_{w_2}V)}{N_0^T(\pa_{w_2}V\times V)}-B_1\frac{N_0^T(\pa_{w_1}V\times V)^2}{N_0^T(\pa_{w_2}V\times V)^2}\frac{N_0^T(\pa_{w_2}^2V\times\pa_{w_2}V)}{N_0^T(\pa_{w_2}V\times V)}-2B_1\frac{N_0^T(\pa_{w_1}V\times V)}{N_0^T(\pa_{w_2}V\times V)}\frac{N_0^T(\pa_{w_1w_2}^2V\times\pa_{w_2}V)}{N_0^T(\pa_{w_2}V\times V)})]$$
$$+\frac{\pa_{w_1}\mathbf{m}}{\mathcal{C}_2}[-\frac{\pa_{w_2}\mathbf{m}}{\mathcal{C}_2}\frac{N_0^T[V\times d(V+x_0)]}{N_0^T(\pa_{w_2}V\times V)}\wedge\frac{N_0^T(\pa_{w_2}V\times dV)}{N_0^T(\pa_{w_2}V\times V)}
+\frac{N_0^T(V\times dV)}{N_0^T(\pa_{w_2}V\times V)}\wedge$$$$ B_1\frac{N_0^T(\pa_{w_1}V\times \pa_{w_2}V)N_0^T(V\times dx_0)}{N_0^T(\pa_{w_2}V\times V)\mathcal{C}_2}[-\frac{N_0^T(\pa_{w_1}V\times \pa_{w_2}V)}{N_0^T(\pa_{w_2}V\times V)}\frac{\pa_{w_2}\mathbf{m}}{\mathcal{C}_2}+\frac{\pa_{w_2}\mathbf{m}}{\mathcal{C}_2}
\frac{N_0^T(\pa_{w_1w_2}^2V\times V)}{N_0^T(\pa_{w_2}V\times V)}$$$$-\frac{N_0^T(\pa_{w_2}^2V\times V)}{N_0^T(\pa_{w_2}V\times V)}(\mathbf{m}\frac{N_0^T(\pa_{w_1}V\times\pa_{w_2}V)}{N_0^T(\pa_{w_2}V\times V)}+\pa_{w_1}\mathbf{m})]+$$$$\frac{N_0^T(\pa_{w_2}V\times dV)}{N_0^T(\pa_{w_2}V\times V)}\wedge B_1\frac{N_0^T(\pa_{w_1}V\times \pa_{w_2}V)N_0^T(V\times dx_0)}{N_0^T(\pa_{w_2}V\times V)\mathcal{C}_2^2}[-\mathbf{m}\frac{N_0^T(\pa_{w_1}V\times\pa_{w_2}V)}{N_0^T(\pa_{w_2}V\times V)}+$$$$2\pa_{w_2}\mathbf{m}\frac{N_0^T(\pa_{w_1}V\times V)}{N_0^T(\pa_{w_2}V\times V)}-2\pa_{w_1}\mathbf{m}-\mathbf{m}\frac{N_0^T(\pa_{w_1w_2}^2V\times V)}{N_0^T(\pa_{w_2}V\times V)}+\mathbf{m}\frac{N_0^T(\pa_{w_2}^2V\times V)}{N_0^T(\pa_{w_2}V\times V)}\frac{N_0^T(\pa_{w_1}V\times V)}{N_0^T(\pa_{w_2}V\times V)}]$$$$-\frac{N_0^T(d\pa_{w_2}V\times V)}{N_0^T(\pa_{w_2}V\times V)}\wedge
B_1\frac{N_0^T(\pa_{w_1}V\times \pa_{w_2}V)N_0^T(V\times dx_0)}{N_0^T(\pa_{w_2}V\times V)\mathcal{C}_2}]+\frac{\pa_{w_2}\mathbf{m}}{\mathcal{C}_2}[\frac{\mathbf{m}}{\mathcal{C}_2}
\frac{N_0^T(\pa_{w_2}V\times d\pa_{w_1}V)}{N_0^T(\pa_{w_2}V\times V)}\wedge$$$$\frac{N_0^T[V\times d(V+x_0)]}{N_0^T(\pa_{w_2}V\times V)}-\frac{\mathbf{m}}{\mathcal{C}_2}
\frac{N_0^T(\pa_{w_2}V\times d\pa_{w_2}V)}{N_0^T(\pa_{w_2}V\times V)}\wedge\frac{N_0^T[V\times d(V+x_0)]}{N_0^T(\pa_{w_2}V\times V)}-\frac{N_0^T(V\times d\pa_{w_2}V)}{N_0^T(\pa_{w_2}V\times V)}\wedge
$$$$\frac{N_0^T(\pa_{w_1}V\times V)}{N_0^T(\pa_{w_2}V\times V)}B_1\frac{N_0^T(\pa_{w_1}V\times \pa_{w_2}V)N_0^T(V\times dx_0)}{N_0^T(\pa_{w_2}V\times V)\mathcal{C}_2}+\frac{N_0^T[V\times d(V+x_0)]}{N_0^T(\pa_{w_2}V\times V)}\wedge$$$$[B_1\frac{N_0^T(\pa_{w_1}V\times \pa_{w_2}V)^2N_0^T(V\times dx_0)}{N_0^T(\pa_{w_2}V\times V)^2\mathcal{C}_2}-\pa_{w_1}B_1\frac{N_0^T(\pa_{w_1}V\times \pa_{w_2}V)N_0^T(V\times dx_0)}{N_0^T(\pa_{w_2}V\times V)\mathcal{C}_2}-$$$$B_1\frac{N_0^T(\pa_{w_1}^2V\times \pa_{w_2}V)N_0^T(V\times dx_0)}{N_0^T(\pa_{w_2}V\times V)\mathcal{C}_2}-B_1\frac{N_0^T(\pa_{w_1}V\times \pa_{w_1w_2}^2V)N_0^T(V\times dx_0)}{N_0^T(\pa_{w_2}V\times V)\mathcal{C}_2}+$$$$\frac{\mathbf{m}}{\mathcal{C}_2}\frac{N_0^T(\pa_{w_1w_2}^2V\times \pa_{w_2}V)}{N_0^T(\pa_{w_2}V\times V)}\frac{N_0^T(V\times dV)}{N_0^T(\pa_{w_2}V\times V)}+\frac{N_0^T(\pa_{w_2}^2V\times V)}{N_0^T(\pa_{w_2}V\times V)}\frac{N_0^T(\pa_{w_1}V\times V)}{N_0^T(\pa_{w_2}V\times V)}B_1$$$$\frac{N_0^T(\pa_{w_1}V\times \pa_{w_2}V)N_0^T(V\times dx_0)}{N_0^T(\pa_{w_2}V\times V)\mathcal{C}_2}-\frac{N_0^T(\pa_{w_2}V\times dx_0)}{N_0^T(\pa_{w_2}V\times V)\mathcal{C}_2}(\mathbf{m}\frac{N_0^T(\pa_{w_1}V\times\pa_{w_2}V)}{N_0^T(\pa_{w_2}V\times V)}+\pa_{w_1}\mathbf{m})$$$$-\frac{\pa_{w_2}\mathbf{m}}{\mathcal{C}_2}\frac{N_0^T(\pa_{w_1}V\times V)}{N_0^T(\pa_{w_2}V\times V)}B_1\frac{N_0^T(\pa_{w_1}V\times \pa_{w_2}V)N_0^T(V\times dx_0)}{N_0^T(\pa_{w_2}V\times V)\mathcal{C}_2}(\frac{N_0^T(\pa_{w_1w_2}^2V\times V)}{N_0^T(\pa_{w_2}V\times V)}-$$$$\frac{N_0^T(\pa_{w_1}V\times V)}{N_0^T(\pa_{w_2}V\times V)}\frac{N_0^T(\pa_{w_2}^2V\times V)}{N_0^T(\pa_{w_2}V\times V)})+\frac{N_0^T(\pa_{w_1}V\times V)}{N_0^T(\pa_{w_2}V\times V)}(\pa_{w_2}B_1\frac{N_0^T(\pa_{w_1}V\times \pa_{w_2}V)N_0^T(V\times dx_0)}{N_0^T(\pa_{w_2}V\times V)\mathcal{C}_2}+$$$$B_1\frac{N_0^T(\pa_{w_1w_2}^2V\times \pa_{w_2}V)N_0^T(V\times dx_0)}{N_0^T(\pa_{w_2}V\times V)\mathcal{C}_2}+B_1\frac{N_0^T(\pa_{w_1}V\times \pa_{w_2}^2V)N_0^T(V\times dx_0)}{N_0^T(\pa_{w_2}V\times V)\mathcal{C}_2})-$$$$\frac{\mathbf{m}}{\mathcal{C}_2}\frac{N_0^T(\pa_{w_2}^2V\times\pa_{w_2}V)}{N_0^T(\pa_{w_2}V\times V)}\frac{N_0^T(V\times dV)}{N_0^T(\pa_{w_2}V\times V)}-\frac{N_0^T(\pa_{w_2}V\times dV)}{N_0^T(\pa_{w_2}V\times V)\mathcal{C}_2}(2\mathbf{m}\frac{N_0^T(\pa_{w_1}V\times\pa_{w_2}V)}{N_0^T(\pa_{w_2}V\times V)}+\pa_{w_1}\mathbf{m})$$$$+\frac{\pa_{w_2}\mathbf{m}}{\mathcal{C}_2}\frac{N_0^T[\pa_{w_1}V\times d(2V+x_0)]}{N_0^T(\pa_{w_2}V\times V)}-(\frac{N_0^T(\pa_{w_1w_2}^2V\times V)}{N_0^T(\pa_{w_2}V\times V)}-\frac{N_0^T(\pa_{w_1}V\times V)}{N_0^T(\pa_{w_2}V\times V)}\frac{N_0^T(\pa_{w_2}^2V\times V)}{N_0^T(\pa_{w_2}V\times V)})\frac{\pa_{w_2}\mathbf{m}}{\mathcal{C}_2}$$$$\frac{N_0^T(\pa_{w_1}V\times V)}{N_0^T(\pa_{w_2}V\times V)}B_1\frac{N_0^T(\pa_{w_1}V\times \pa_{w_2}V)N_0^T(V\times dx_0)}{N_0^T(\pa_{w_2}V\times V)\mathcal{C}_2}+\frac{N_0^T(\pa_{w_1}V\times\pa_{w_2}V)}{N_0^T(\pa_{w_2}V\times V)}\frac{1}{\mathcal{C}_2}(\mathbf{m}\frac{N_0^T(\pa_{w_1}V\times\pa_{w_2}V)}{N_0^T(\pa_{w_2}V\times V)}$$$$+\pa_{w_1}\mathbf{m})B_1\frac{N_0^T(\pa_{w_1}V\times \pa_{w_2}V)N_0^T(V\times dx_0)}{N_0^T(\pa_{w_2}V\times V)\mathcal{C}_2}]+\frac{N_0^T(\pa_{w_2}V\times dV)}{N_0^T(\pa_{w_2}V\times V)}\wedge[-\frac{\mathbf{m}}{\mathcal{C}_2}\frac{N_0^T(\pa_{w_1}V\times dV)}{N_0^T(\pa_{w_2}V\times V)}+\frac{N_0^T(\pa_{w_1}V\times V)}{N_0^T(\pa_{w_2}V\times V)}$$$$B_1\frac{N_0^T(\pa_{w_1}V\times \pa_{w_2}V)N_0^T(V\times dx_0)}{N_0^T(\pa_{w_2}V\times V)\mathcal{C}_2}[1+\frac{\mathbf{m}}{\mathcal{C}_2}(\frac{N_0^T(\pa_{w_1w_2}^2V\times V)}{N_0^T(\pa_{w_2}V\times V)}-\frac{N_0^T(\pa_{w_1}V\times V)}{N_0^T(\pa_{w_2}V\times V)}\frac{N_0^T(\pa_{w_2}^2V\times V)}{N_0^T(\pa_{w_2}V\times V)}-$$$$\frac{N_0^T(\pa_{w_1}V\times\pa_{w_2}V)}{N_0^T(\pa_{w_2}V\times V)})]]+\frac{N_0^T[\pa_{w_1}V\times d(2V+x_0)]}{N_0^T(\pa_{w_2}V\times V)}\wedge B_1\frac{N_0^T(\pa_{w_1}V\times \pa_{w_2}V)N_0^T(V\times dx_0)}{N_0^T(\pa_{w_2}V\times V)\mathcal{C}_2}-$$$$\frac{N_0^T[dV\times\wedge d(V+x_0)]}{N_0^T(\pa_{w_2}V\times V)}+\frac{1}{N_0^T(\pa_{w_2}V\times V)}|V|^2\frac{1}{2}N_0^T(dN_0\times\wedge dN_0)]+\frac{\mathbf{m}}{\mathcal{C}_2}[\frac{N_0^T(V\times d\pa_{w_2}V)}{N_0^T(\pa_{w_2}V\times V)}\wedge $$$$B_1\frac{N_0^T(\pa_{w_1}V\times \pa_{w_2}V)^2N_0^T(V\times dx_0)}{N_0^T(\pa_{w_2}V\times V)^2\mathcal{C}_2}-\frac{N_0^T(\pa_{w_2}V\times d\pa_{w_2}V)}{N_0^T(\pa_{w_2}V\times V)}\wedge\frac{N_0^T(V\times dx_0)}{N_0^T(\pa_{w_2}V\times V)}+\frac{N_0^T[V\times d(V+x_0)]}{N_0^T(\pa_{w_2}V\times V)}$$$$\wedge[-\frac{N_0^T(\pa_{w_2}^2V\times\pa_{w_2}V)}{N_0^T(\pa_{w_2}V\times V)}\frac{N_0^T(V\times dV)}{N_0^T(\pa_{w_2}V\times V)}-\frac{\mathbf{m}}{\mathcal{C}_2}\frac{N_0^T(\pa_{w_2}V\times dV)}{N_0^T(\pa_{w_2}V\times V)}(\frac{N_0^T(\pa_{w_1w_2}^2V\times\pa_{w_2}V)}{N_0^T(\pa_{w_2}V\times V)}-$$$$\frac{N_0^T(\pa_{w_2}^2V\times\pa_{w_2}V)}{N_0^T(\pa_{w_2}V\times V)}\frac{N_0^T(\pa_{w_1}V\times V)}{N_0^T(\pa_{w_2}V\times V)})-B_1\frac{N_0^T(\pa_{w_1}V\times \pa_{w_2}V)N_0^T(V\times dx_0)}{N_0^T(\pa_{w_2}V\times V)^2\mathcal{C}_2}(N_0^T(\pa_{w_1w_2}^2V\times\pa_{w_2}V)+$$$$N_0^T(\pa_{w_1}V\times\pa_{w_2}^2V))]+
\frac{N_0^T(\pa_{w_2}V\times dV)}{N_0^T(\pa_{w_2}V\times V)}\wedge[[\frac{N_0^T(\pa_{w_1w_2}^2V\times V)}{N_0^T(\pa_{w_2}V\times V)}-\frac{N_0^T(\pa_{w_1}V\times\pa_{w_2}V)}{N_0^T(\pa_{w_2}V\times V)}+$$$$\frac{\mathbf{m}}{\mathcal{C}_2}(\frac{N_0^T(\pa_{w_1}^2V\times\pa_{w_2}V)}{N_0^T(\pa_{w_2}V\times V)}+\frac{N_0^T(\pa_{w_1}V\times\pa_{w_2}V)^2}{N_0^T(\pa_{w_2}V\times V)^2})]B_1\frac{N_0^T(\pa_{w_1}V\times \pa_{w_2}V)N_0^T(V\times dx_0)}{N_0^T(\pa_{w_2}V\times V)\mathcal{C}_2}+$$$$\frac{\mathbf{m}}{\mathcal{C}_2}[\frac{N_0^T(\pa_{w_1w_2}^2V\times dV)}{N_0^T(\pa_{w_2}V\times V)}-B_1\frac{N_0^T(\pa_{w_1}V\times\pa_{w_2}V)N_0^T(V\times dx_0)}{N_0^T(\pa_{w_2}V\times V)^2\mathcal{C}_2}(N_0^T(\pa_{w_1}^2V\times\pa_{w_2}V)+$$$$
N_0^T(\pa_{w_1}V\times\pa_{w_1w_2}^2V))]+
\frac{1}{\mathcal{C}_2}(\pa_{w_1}B_1\frac{N_0^T(\pa_{w_1}V\times\pa_{w_2}V)N_0^T(V\times dx_0)}{N_0^T(\pa_{w_2}V\times V)}+$$$$B_1\frac{N_0^T(\pa_{w_1}^2V\times\pa_{w_2}V)N_0^T(V\times dx_0)}{N_0^T(\pa_{w_2}V\times V)}+B_1\frac{N_0^T(\pa_{w_1}V\times\pa_{w_1w_2}^2V)N_0^T(V\times dx_0)}{N_0^T(\pa_{w_2}V\times V)}-$$$$B_1\frac{N_0^T(\pa_{w_1}V\times\pa_{w_2}V)^2N_0^T(V\times dx_0)}{N_0^T(\pa_{w_2}V\times V)^2})-(\frac{N_0^T(\pa_{w_2}^2V\times V)}{N_0^T(\pa_{w_2}V\times V)}-\frac{\mathbf{m}}{\mathcal{C}_2}\frac{N_0^T(\pa_{w_1}V\times V)}{N_0^T(\pa_{w_2}V\times V)}\frac{N_0^T(\pa_{w_2}^2V\times\pa_{w_2}V)}{N_0^T(\pa_{w_2}V\times V)})$$$$\frac{N_0^T(\pa_{w_1}V\times V)}{N_0^T(\pa_{w_2}V\times V)}B_1\frac{N_0^T(\pa_{w_1}V\times\pa_{w_2}V)N_0^T(V\times dx_0)}{N_0^T(\pa_{w_2}V\times V)\mathcal{C}_2}+\frac{\mathbf{m}}{\mathcal{C}_2}\frac{N_0^T(\pa_{w_1}V\times V)}{N_0^T(\pa_{w_2}V\times V)}[-\frac{N_0^T(\pa_{w_2}^2V\times dV)}{N_0^T(\pa_{w_2}V\times V)}+$$$$B_1\frac{N_0^T(\pa_{w_1}V\times\pa_{w_2}V)N_0^T(V\times dx_0)}{N_0^T(\pa_{w_2}V\times V)^2\mathcal{C}_2}(N_0^T(\pa_{w_1w_2}^2V\times\pa_{w_2}V)+N_0^T(\pa_{w_1}V\times\pa_{w_2}^2V))]
-\frac{1}{\mathcal{C}_2}\frac{N_0^T(\pa_{w_1}V\times V)}{N_0^T(\pa_{w_2}V\times V)}$$$$(\pa_{w_2}B_1\frac{N_0^T(\pa_{w_1}V\times\pa_{w_2}V)N_0^T(V\times dx_0)}{N_0^T(\pa_{w_2}V\times V)}
+B_1\frac{N_0^T(\pa_{w_1w_2}^2V\times\pa_{w_2}V)N_0^T(V\times dx_0)}{N_0^T(\pa_{w_2}V\times V)}+$$$$
B_1\frac{N_0^T(\pa_{w_1}V\times\pa_{w_2}^2V)N_0^T(V\times dx_0)}{N_0^T(\pa_{w_2}V\times V)})]
-\frac{N_0^T(\pa_{w_1w_2}^2V\times dV)}{N_0^T(\pa_{w_2}V\times V)}\wedge B_1\frac{N_0^T(\pa_{w_1}V\times\pa_{w_2}V)N_0^T(V\times dx_0)}{N_0^T(\pa_{w_2}V\times V)\mathcal{C}_2}$$$$
+\frac{N_0^T(\pa_{w_2}^2V\times dV)}{N_0^T(\pa_{w_2}V\times V)}\wedge\frac{N_0^T(\pa_{w_1}V\times V)}{N_0^T(\pa_{w_2}V\times V)}B_1\frac{N_0^T(\pa_{w_1}V\times\pa_{w_2}V)N_0^T(V\times dx_0)}{N_0^T(\pa_{w_2}V\times V)\mathcal{C}_2}-\frac{1}{2}\frac{\pa_{w_2}V^TV}{N_0^T(\pa_{w_2}V\times V)}$$$$N_0^T(dN_0\times\wedge dN_0)]-(\pa_{w_2}B_1\frac{N_0^T(\pa_{w_1}V\times\pa_{w_2}V)N_0^T(V\times dx_0)}{N_0^T(\pa_{w_2}V\times V)\mathcal{C}_2}+B_1\frac{N_0^T(\pa_{w_1w_2}^2V\times\pa_{w_2}V)N_0^T(V\times dx_0)}{N_0^T(\pa_{w_2}V\times V)\mathcal{C}_2}$$$$+B_1\frac{N_0^T(\pa_{w_1}V\times\pa_{w_2}^2V)N_0^T(V\times dx_0)}{N_0^T(\pa_{w_2}V\times V)\mathcal{C}_2})\wedge\frac{N_0^T(V\times dV)}{N_0^T(\pa_{w_2}V\times V)}-B_1\frac{N_0^T(\pa_{w_1}V\times\pa_{w_2}V)N_0^T(\pa_{w_2}V\times dx_0)}{N_0^T(\pa_{w_2}V\times V)\mathcal{C}_2}$$$$\wedge\frac{N_0^T[V\times d(V+x_0)]}{N_0^T(\pa_{w_2}V\times V)}-\frac{1}{2}\mathbf{m}B_1N_0^T(dN_0\times\wedge dN_0)+dB_1\wedge\frac{N_0^T(\pa_{w_1}V\times\pa_{w_2}V)N_0^T(V\times dx_0)}{N_0^T(\pa_{w_2}V\times V)\mathcal{C}_2}+$$$$B_1\frac{N_0^T(d\pa_{w_1}V\times\pa_{w_2}V)\wedge N_0^T(V\times dx_0)}{N_0^T(\pa_{w_2}V\times V)\mathcal{C}_2}+B_1\frac{N_0^T(\pa_{w_1}V\times d\pa_{w_2}V)\wedge N_0^T(V\times dx_0)}{N_0^T(\pa_{w_2}V\times V)\mathcal{C}_2}+$$$$B_1\frac{N_0^T(\pa_{w_1}V\times\pa_{w_2}V)N_0^T(dV\times\wedge dx_0)}{N_0^T(\pa_{w_2}V\times V)\mathcal{C}_2}+
[\frac{d\pa_{w_1}\mathbf{m}}{\mathcal{C}_2}B_1+\frac{d\pa_{w_2}\mathbf{m}}{\mathcal{C}_2}
(\frac{\mathbf{m}}{N_0^T(\pa_{w_2}V\times V)}-B_1\frac{N_0^T(\pa_{w_1}V\times V)}{N_0^T(\pa_{w_2}V\times V)})-$$$$\frac{\pa_{w_1}^2\mathbf{m}}{\mathcal{C}_2}\frac{\mathcal{B}_2}{\mathcal{C}_2}B_1+
\frac{\pa_{w_2}^2\mathbf{m}}{\mathcal{C}_2}[\frac{N_0^T[V\times d(V+x_0)]}{N_0^T(\pa_{w_2}V\times V)}+\frac{N_0^T(\pa_{w_1}V\times V)}{N_0^T(\pa_{w_2}V\times V)}(\frac{\mathcal{B}_2}{\mathcal{C}_2}-B_1\frac{N_0^T(\pa_{w_1}V\times \pa_{w_2}V)N_0^T(V\times dx_0)}
{N_0^T(\pa_{w_2}V\times V)\mathcal{C}_2})]$$
$$(\frac{\mathbf{m}}{N_0^T(\pa_{w_2}V\times V)}-B_1\frac{N_0^T(\pa_{w_1}V\times V)}{N_0^T(\pa_{w_2}V\times V)})+\frac{\pa_{w_1w_2}^2\mathbf{m}}{\mathcal{C}_2}[(\frac{N_0^T[V\times d(V+x_0)]}{N_0^T(\pa_{w_2}V\times V)}+\frac{N_0^T(\pa_{w_1}V\times V)}{N_0^T(\pa_{w_2}V\times V)}(\frac{\mathcal{B}_2}{\mathcal{C}_2}-$$$$B_1\frac{N_0^T(\pa_{w_1}V\times \pa_{w_2}V)N_0^T(V\times dx_0)}
{N_0^T(\pa_{w_2}V\times V)\mathcal{C}_2}))B_1-(\frac{\mathcal{B}_2}{\mathcal{C}_2}-B_1\frac{N_0^T(\pa_{w_1}V\times \pa_{w_2}V)N_0^T(V\times dx_0)}
{N_0^T(\pa_{w_2}V\times V)\mathcal{C}_2})$$$$(\frac{\mathbf{m}}{N_0^T(\pa_{w_2}V\times V)}-B_1\frac{N_0^T(\pa_{w_1}V\times V)}{N_0^T(\pa_{w_2}V\times V)})]+\frac{d\mathbf{m}}{\mathcal{C}_2}[-\frac{B_1}{\mathcal{C}_2}(\pa_{w_1}\mathbf{m}
\frac{N_0^T(\pa_{w_1}V\times\pa_{w_2}V)}{N_0^T(\pa_{w_2}V\times V)}+\mathbf{m}\frac{N_0^T(\pa_{w_1}^2V\times\pa_{w_2}V)}{N_0^T(\pa_{w_2}V\times V)}$$$$-2B_1\frac{\mathbf{m}}{\mathcal{C}_2}\frac{N_0^T(\pa_{w_1}V\times V)}{N_0^T(\pa_{w_2}V\times V)}\frac{N_0^T(\pa_{w_1w_2}^2V\times\pa_{w_2}V)}{N_0^T(\pa_{w_2}V\times V)}+\mathbf{m}\frac{N_0^T(\pa_{w_1}V\times\pa_{w_2}V)^2}{N_0^T(\pa_{w_2}V\times V)^2}-\pa_{w_2}\mathbf{m}\frac{N_0^T(\pa_{w_1}^2V\times V)}{N_0^T(\pa_{w_2}V\times V)}+$$$$2\pa_{w_2}\mathbf{m}\frac{N_0^T(\pa_{w_1}V\times V)}{N_0^T(\pa_{w_2}V\times V)}\frac{N_0^T(\pa_{w_1w_2}^2V\times V)}{N_0^T(\pa_{w_2}V\times V)}-2\pa_{w_2}\mathbf{m}\frac{N_0^T(\pa_{w_1}V\times V)}{N_0^T(\pa_{w_2}V\times V)}\frac{N_0^T(\pa_{w_1}V\times \pa_{w_2}V)}{N_0^T(\pa_{w_2}V\times V)})-$$$$\pa_{w_1}B_1-\frac{\mathbf{m}}{N_0^T(\pa_{w_2}V\times V)\mathcal{C}_2}(\pa_{w_2}\mathbf{m}\frac{N_0^T(\pa_{w_1}V\times\pa_{w_2}V)}{N_0^T(\pa_{w_2}V\times V)}+\mathbf{m}\frac{N_0^T(\pa_{w_1w_2}^2V\times\pa_{w_2}V)}{N_0^T(\pa_{w_2}V\times V)}+\mathbf{m}\frac{N_0^T(\pa_{w_1}V\times\pa_{w_2}^2V)}{N_0^T(\pa_{w_2}V\times V)}$$$$-\mathbf{m}\frac{N_0^T(\pa_{w_1}V\times\pa_{w_2}V)}{N_0^T(\pa_{w_2}V\times V)}\frac{N_0^T(\pa_{w_2}^2V\times V)}{N_0^T(\pa_{w_2}V\times V)}-\pa_{w_2}\mathbf{m}\frac{N_0^T(\pa_{w_1w_2}^2V\times V)}{N_0^T(\pa_{w_2}V\times V)}-\pa_{w_2}\mathbf{m}\frac{N_0^T(\pa_{w_1}V\times\pa_{w_2}V)}{N_0^T(\pa_{w_2}V\times V)}+$$$$\pa_{w_2}\mathbf{m}\frac{N_0^T(\pa_{w_1}V\times V)}{N_0^T(\pa_{w_2}V\times V)}\frac{N_0^T(\pa_{w_2}^2V\times V)}{N_0^T(\pa_{w_2}V\times V)})+B_1\frac{N_0^T(\pa_{w_1}V\times V)}{N_0^T(\pa_{w_2}V\times V)\mathcal{C}_2}(\mathbf{m}\frac{N_0^T(\pa_{w_1}V\times\pa_{w_2}^2V)}{N_0^T(\pa_{w_2}V\times V)}-$$$$\mathbf{m}\frac{N_0^T(\pa_{w_1}V\times\pa_{w_2}V)}{N_0^T(\pa_{w_2}V\times V)}\frac{N_0^T(\pa_{w_2}^2V\times V)}{N_0^T(\pa_{w_2}V\times V)}+\pa_{w_2}\mathbf{m}\frac{N_0^T(\pa_{w_1}V\times V)}{N_0^T(\pa_{w_2}V\times V)}\frac{N_0^T(\pa_{w_2}^2V\times V)}{N_0^T(\pa_{w_2}V\times V)})+\pa_{w_2}B_1\frac{N_0^T(\pa_{w_1}V\times V)}{N_0^T(\pa_{w_2}V\times V)}]$$$$+\frac{\pa_{w_1}\mathbf{m}}{\mathcal{C}_2}[-B_1\frac{N_0^T(\pa_{w_2}V\times dV)}{N_0^T(\pa_{w_2}V\times V)}-B_1\frac{N_0^T(\pa_{w_1}V\times\pa_{w_2}V)}{N_0^T(\pa_{w_2}V\times V)}\frac{\pa_{w_2}\mathbf{m}\frac{N_0^T[V\times d(V+x_0)]}{N_0^T(\pa_{w_2}V\times V)}-
\mathbf{m}\frac{N_0^T(\pa_{w_2}V\times dV)}{N_0^T(\pa_{w_2}V\times V)}}{\mathcal{C}_2}$$$$+B_1^2\frac{N_0^T(\pa_{w_1}V\times \pa_{w_2}V)N_0^T(V\times dx_0)}
{N_0^T(\pa_{w_2}V\times V)^2\mathcal{C}_2}N_0^T(\pa_{w_1w_2}^2V\times V)+B_1\frac{N_0^T(\pa_{w_1}V\times \pa_{w_2}V)N_0^T(V\times dx_0)}
{N_0^T(\pa_{w_2}V\times V)^2\mathcal{C}_2}$$$$N_0^T(\pa_{w_2}^2V\times V)(\frac{\mathbf{m}}{N_0^T(\pa_{w_2}V\times V)}-B_1\frac{N_0^T(\pa_
{w_1}V\times V)}{N_0^T(\pa_{w_2}V\times V)})]+\frac{\pa_{w_2}\mathbf{m}}{\mathcal{C}_2}[B_1\frac{N_0^T(V\times d\pa_{w_1}V)}{N_0^T(\pa_{w_2}V\times V)}+$$$$\frac{N_0^T(V\times d\pa_{w_2}V)}{N_0^T(\pa_{w_2}V\times V)}
(\frac{\mathbf{m}}{N_0^T(\pa_{w_2}V\times V)}-B_1\frac{N_0^T(\pa_{w_1}V\times V)}
{N_0^T(\pa_{w_2}V\times V)})+\frac{N_0^T[\pa_{w_2}V\times d(V+x_0)]}{N_0^T(\pa_{w_2}V\times V)}\frac{\mathbf{m}}{N_0^T(\pa_{w_2}V\times V)}$$$$+\frac{N_0^T[V\times d(V+x_0)]}{N_0^T(\pa_{w_2}V\times V)}[-B_1\frac{N_0^T(\pa_{w_1w_2}^2V\times V)}{N_0^T(\pa_{w_2}V\times V)}+\frac{\pa_{w_2}\mathbf{m}}{\mathcal{C}_2}(B_1\frac{N_0^T(\pa_{w_1}^2V\times V)}{N_0^T(\pa_{w_2}V\times V)}-B_1\frac{N_0^T(\pa_{w_1}V\times V)}{N_0^T(\pa_{w_2}V\times V)}$$$$\frac{N_0^T(\pa_{w_1w_2}^2V\times V)}{N_0^T(\pa_{w_2}V\times V)}+B_1\frac{N_0^T(\pa_{w_1}V\times V)}{N_0^T(\pa_{w_2}V\times V)}\frac{N_0^T(\pa_{w_1}V\times\pa_{w_2}V)}{N_0^T(\pa_{w_2}V\times V)}+(\frac{N_0^T(\pa_{w_1w_2}^2V\times V)}{N_0^T(\pa_{w_2}V\times V)}-\frac{N_0^T(\pa_{w_1}V\times V)}{N_0^T(\pa_{w_2}V\times V)}$$$$\frac{N_0^T(\pa_{w_2}^2V\times V)}{N_0^T(\pa_{w_2}V\times V)})(\frac{\mathbf{m}}{N_0^T(\pa_{w_2}V\times V)}-B_1\frac{N_0^T(\pa_{w_1}V\times V)}{N_0^T(\pa_{w_2}V\times V)}))-\pa_{w_1}B_1-\frac{N_0^T(\pa_{w_2}^2V\times V)}{N_0^T(\pa_{w_2}V\times V)}$$$$(\frac{\mathbf{m}}{N_0^T(\pa_{w_2}V\times V)}-B_1\frac{N_0^T(\pa_{w_1}V\times V)}{N_0^T(\pa_{w_2}V\times V)})+\pa_{w_2}B_1\frac{N_0^T(\pa_{w_1}V\times V)}{N_0^T(\pa_{w_2}V\times V)}]+\frac{N_0^T(\pa_{w_2}V\times dV)}{N_0^T(\pa_{w_2}V\times V)}$$$$[-\frac{\mathbf{m}}{N_0^T(\pa_{w_2}V\times V)}+B_1\frac{N_0^T(\pa_{w_1}V\times V)}{N_0^T(\pa_{w_2}V\times V)}-\frac{\mathbf{m}}{\mathcal{C}_2}(B_1\frac{N_0^T(\pa_{w_1}^2V\times V)}{N_0^T(\pa_{w_2}V\times V)}-B_1\frac{N_0^T(\pa_{w_1}V\times V)}{N_0^T(\pa_{w_2}V\times V)}$$$$\frac{N_0^T(\pa_{w_1w_2}^2V\times V)}{N_0^T(\pa_{w_2}V\times V)}+B_1\frac{N_0^T(\pa_{w_1}V\times V)}{N_0^T(\pa_{w_2}V\times V)}\frac{N_0^T(\pa_{w_1}V\times\pa_{w_2}V)}{N_0^T(\pa_{w_2}V\times V)}+(\frac{N_0^T(\pa_{w_1w_2}^2V\times V)}{N_0^T(\pa_{w_2}V\times V)}-\frac{N_0^T(\pa_{w_1}V\times V)}{N_0^T(\pa_{w_2}V\times V)}$$$$\frac{N_0^T(\pa_{w_2}^2V\times V)}{N_0^T(\pa_{w_2}V\times V)})(\frac{\mathbf{m}}{N_0^T(\pa_{w_2}V\times V)}-B_1\frac{N_0^T(\pa_{w_1}V\times V)}{N_0^T(\pa_{w_2}V\times V)}))]+B_1\frac{N_0^T(\pa_{w_1}V\times \pa_{w_2}V)N_0^T(V\times dx_0)}{N_0^T(\pa_{w_2}V\times V)\mathcal{C}_2}$$$$[-B_1\frac{N_0^T(\pa_{w_1}^2V\times V)}{N_0^T(\pa_{w_2}V\times V)}-(\frac{N_0^T(\pa_{w_1w_2}^2V\times V)}{N_0^T(\pa_{w_2}V\times V)}-\frac{N_0^T(\pa_{w_1}V\times V)}{N_0^T(\pa_{w_2}V\times V)}\frac{N_0^T(\pa_{w_2}^2V\times V)}{N_0^T(\pa_{w_2}V\times V)})(\frac{\mathbf{m}}{N_0^T(\pa_{w_2}V\times V)}$$$$-B_1\frac{N_0^T(\pa_{w_1}V\times V)}{N_0^T(\pa_{w_2}V\times V)})] ]+\frac{\mathbf{m}}{\mathcal{C}_2}[-B_1\frac{N_0^T(\pa_{w_2}V\times d\pa_{w_1}V)}{N_0^T(\pa_{w_2}V\times V)}-\frac{N_0^T(\pa_{w_2}V\times d\pa_{w_2}V)}{N_0^T(\pa_{w_2}V\times V)}(\frac{\mathbf{m}}{N_0^T(\pa_{w_2}V\times V)}-$$$$B_1\frac{N_0^T(\pa_{w_1}V\times V)}{N_0^T(\pa_{w_2}V\times V)})+\frac{N_0^T(V\times dV)}{N_0^T(\pa_{w_2}V\times V)}[B_1\frac{N_0^T(\pa_{w_1w_2}^2V\times\pa_{w_2}V)}{N_0^T(\pa_{w_2}V\times V)}
+(\frac{\mathbf{m}}{N_0^T(\pa_{w_2}V\times V)}-$$$$B_1\frac{N_0^T(\pa_{w_1}V\times V)}{N_0^T(\pa_{w_2}V\times V)})\frac{N_0^T(\pa_{w_2}^2V\times\pa_{w_2}V)}{N_0^T(\pa_{w_2}V\times V)}]+\frac{N_0^T(\pa_{w_2}V\times dV)}{N_0^T(\pa_{w_2}V\times V)}[-B_1\frac{N_0^T(\pa_{w_1}V\times\pa_{w_2}V)}{N_0^T(\pa_{w_2}V\times V)}+B_1\frac{\mathbf{m}}{\mathcal{C}_2}$$$$(\frac{N_0^T(\pa_{w_1}^2V\times\pa_{w_2}V)}{N_0^T(\pa_{w_2}V\times V)}-\frac{N_0^T(\pa_{w_1}V\times V)}{N_0^T(\pa_{w_2}V\times V)}\frac{N_0^T(\pa_{w_1w_2}^2V\times\pa_{w_2}V)}{N_0^T(\pa_{w_2}V\times V)}+
\frac{N_0^T(\pa_{w_1}V\times\pa_{w_2}V)^2}{N_0^T(\pa_{w_2}V\times V)^2})+B_1\pa_{w_1}B_1$$$$+\frac{\mathbf{m}}{\mathcal{C}_2}(\frac{\mathbf{m}}{N_0^T(\pa_{w_2}V\times V)}-B_1\frac{N_0^T(\pa_{w_1}V\times V)}{N_0^T(\pa_{w_2}V\times V)})(\frac{N_0^T(\pa_{w_1w_2}^2V\times\pa_{w_2}V)}{N_0^T(\pa_{w_2}V\times V)}-\frac{N_0^T(\pa_{w_1}V\times V)}{N_0^T(\pa_{w_2}V\times V)}$$$$\frac{N_0^T(\pa_{w_2}^2V\times\pa_{w_2}V)}{N_0^T(\pa_{w_2}V\times V)})-\pa_{w_2}B_1\frac{N_0^T(\pa_{w_1}V\times V)}{N_0^T(\pa_{w_2}V\times V)}]+B_1\frac{N_0^T(\pa_{w_1}V\times\pa_{w_2}V)N_0^T(V\times dx_0)}{N_0^T(\pa_{w_2}V\times V)\mathcal{C}_2}[$$$$B_1(N_0^T(\pa_{w_1}^2V\times\pa_{w_2}V)+N_0^T(\pa_{w_1}V\times\pa_{w_1w_2}^2V))+
N_0^T(\pa_{w_1w_2}^2V\times\pa_{w_2}V)+N_0^T(\pa_{w_1}V\times\pa_{w_2}^2V)]$$$$
-\pa_{w_2}B_1\frac{N_0^T(\pa_{w_1}V\times\pa_{w_2}V)N_0^T(V\times dx_0)}{N_0^T(\pa_{w_2}V\times V)^2}-
B_1\frac{N_0^T(\pa_{w_1w_2}^2V\times\pa_{w_2}V)N_0^T(V\times dx_0)}{N_0^T(\pa_{w_2}V\times V)^2}-$$$$
B_1\frac{N_0^T(\pa_{w_1}V\times\pa_{w_2}^2V)N_0^T(V\times dx_0)}{N_0^T(\pa_{w_2}V\times V)^2}-
B_1\frac{N_0^T(\pa_{w_1}V\times\pa_{w_2}V)N_0^T(\pa_{w_2}V\times dx_0)}{N_0^T(\pa_{w_2}V\times V)^2}]$$$$
-B_1^2(\frac{N_0^T(\pa_{w_1}^2V\times\pa_{w_2}V)N_0^T(V\times dx_0)}{N_0^T(\pa_{w_2}V\times V)\mathcal{C}_2}+\frac{N_0^T(\pa_{w_1}V\times\pa_{w_1w_2}^2V)N_0^T(V\times dx_0)}{N_0^T(\pa_{w_2}V\times V)\mathcal{C}_2}$$$$+\frac{N_0^T(\pa_{w_1}V\times\pa_{w_2}V)N_0^T(\pa_{w_1}V\times dx_0)}{N_0^T(\pa_{w_2}V\times V)\mathcal{C}_2})+\pa_{w_2}B_1(\frac{N_0^T[V\times d(V+x_0)]}{N_0^T(\pa_{w_2}V\times V)}-$$$$\frac{N_0^T(\pa_{w_1}V\times\pa_{w_2}V)N_0^T(V\times dx_0)}{N_0^T(\pa_{w_2}V\times V)\mathcal{C}_2}\frac{\mathbf{m}}{N_0^T(\pa_{w_2}V\times V)})-B_1(\frac{\mathbf{m}}{N_0^T(\pa_{w_2}V\times V)}-B_1\frac{N_0^T(\pa_{w_1}V\times V)}{N_0^T(\pa_{w_2}V\times V)})$$$$(\frac{N_0^T(\pa_{w_1w_2}^2V\times\pa_{w_2}V)N_0^T(V\times dx_0)}{N_0^T(\pa_{w_2}V\times V)\mathcal{C}_2}+\frac{N_0^T(\pa_{w_1}V\times\pa_{w_2}^2V)N_0^T(V\times dx_0)}{N_0^T(\pa_{w_2}V\times V)\mathcal{C}_2}+$$$$\frac{N_0^T(\pa_{w_1}V\times\pa_{w_2}V)N_0^T(\pa_{w_2}V\times dx_0)}{N_0^T(\pa_{w_2}V\times V)\mathcal{C}_2})+dB_1+B_1\frac{N_0^T[\pa_{w_2}V\times d(V+x_0)]}{N_0^T(\pa_{w_2}V\times V)}-$$$$B_1(-\frac{\mathcal{B}_2}{\mathcal{C}_2}+B_1\frac{N_0^T(\pa_{w_1}V\times \pa_{w_2}V)N_0^T(V\times dx_0)}{N_0^T(\pa_{w_2}V\times V)\mathcal{C}_2})\frac{N_0^T(\pa_{w_1}V\times \pa_{w_2}V)}{N_0^T(\pa_{w_2}V\times V)}]\wedge V^T(\om\times N_0+dN_0)$$

for all $\om$ satisfying (\ref{eq:om}).

The $4$-dimensional generic integrable rolling distribution of contact elements above splits into an $1$-dimensional family
of generic $3$-dimensional integrable rolling distributions of contact elements
if and only if along the leaves we have
$0=dw_1\wedge dw_2=dw_1\wedge(\frac{N_0^T[V\times d(V+x_0)]}{N_0^T(\pa_{w_2}V\times V)}+
\mathbf{m}\frac{V^T(\om\times N_0+dN_0)}{N_0^T(\pa_{w_2}V\times V)})$, so along the leaves we need
$dw_1=\frac{B_1}{\mathbf{m}}[N_0^T[V\times d(V+x_0)]+\mathbf{m}V^T(\om\times N_0+dN_0)]$.

From (\ref{eq:dw1}) we get $A_1=\frac{B_1}{\mathbf{m}}N_0^T[V\times d(V+x_0)]$ and from (\ref{eq:a1}) we get

$$\frac{\mathcal{B}_2}{\mathcal{C}_2}=(-\frac{N_0^T[V\times d(V+x_0)]}{\mathbf{m}}+
\frac{N_0^T(\pa_{w_1}V\times \pa_{w_2}V)N_0^T(V\times dx_0)}{N_0^T(\pa_{w_2}V\times V)\mathcal{C}_2})B_1,$$

so

\begin{eqnarray}\label{eq:B1}
\mathcal{B}_2\wedge(-\frac{N_0^T[V\times d(V+x_0)]}{\mathbf{m}}+
\frac{N_0^T(\pa_{w_1}V\times \pa_{w_2}V)N_0^T(V\times dx_0)}{N_0^T(\pa_{w_2}V\times V)\mathcal{C}_2})=0,\nonumber\\
B_1=\frac{\pa_u\mathbf{m}+\pa_{w_2}\mathbf{m}\frac{N_0^T[V\times\pa_u(V+x_0)]}{N_0^T(\pa_{w_2}V\times V)}-
\mathbf{m}\frac{N_0^T(\pa_{w_2}V\times\pa_uV)}{N_0^T(\pa_{w_2}V\times V)}}{-\mathcal{C}_2\frac{N_0^T[V\times\pa_u(V+x_0)]}{\mathbf{m}}+
\frac{N_0^T(\pa_{w_1}V\times\pa_{w_2}V)N_0^T(V\times\pa_ux_0)}{N_0^T(\pa_{w_2}V\times V)}}.
\end{eqnarray}

If we exchange the role played by $w_1$ and $w_2$ in the first equation of (\ref{eq:B1}), then the new equation
remains the same as the old one.

The question appears if the contact elements tangent planes to a focal surface of a generic (not normal) Weingarten congruence can be uniquely extended to a $3$-dimensional integrable rolling distribution of contact elements with the symmetry of the tangency configuration.

\end{document}